\documentclass[a4paper, 11pt]{article}

\usepackage{graphicx}
\usepackage{subfig}
\usepackage{color}
\usepackage{bm}
\usepackage{mathrsfs}
\usepackage{ulem}
\usepackage{multirow}
\usepackage{pdflscape}
\usepackage{rotating}
\usepackage{lscape}
\usepackage{pdflscape}
\usepackage{titlesec}

\usepackage{hyperref}

\titleclass{\subsubsubsection}{straight}[\subsection]

\newcounter{subsubsubsection}[subsubsection]
\renewcommand\thesubsubsubsection{\thesubsubsection.\arabic{subsubsubsection}}

\titleformat{\subsubsubsection}
  {\normalfont\normalsize\bfseries}{\thesubsubsubsection}{1em}{}
\titlespacing*{\subsubsubsection}
{0pt}{3.25ex plus 1ex minus .2ex}{1.5ex plus .2ex}

\makeatletter
\renewcommand\paragraph{\@startsection{paragraph}{5}{\z@}%
  {3.25ex \@plus1ex \@minus.2ex}%
  {-1em}%
  {\normalfont\normalsize\bfseries}}
\renewcommand\subparagraph{\@startsection{subparagraph}{6}{\parindent}%
  {3.25ex \@plus1ex \@minus .2ex}%
  {-1em}%
  {\normalfont\normalsize\bfseries}}
\def\toclevel@subsubsubsection{4}
\def\toclevel@paragraph{5}
\def\toclevel@paragraph{6}
\def\l@subsubsubsection{\@dottedtocline{4}{7em}{4em}}
\def\l@paragraph{\@dottedtocline{5}{10em}{5em}}
\def\l@subparagraph{\@dottedtocline{6}{14em}{6em}}
\makeatother

\def\b{\begin{eqnarray}}
\def\e{\end{eqnarray}}
\def\n{\noindent}
\def\p{p_n(x)}
\def\q{q_n(x)}
\def\r{r_{n-1}(x)}

\def\eq{p_n(x) = 0}

\makeatletter
\newcommand*\bigdot{\mathpalette\bigdot@{.5}}
\newcommand*\bigdot@[2]{\mathbin{\vcenter{\hbox{\scalebox{#2}{$\m@th#1\bullet$}}}}}
\makeatother

\begin{document}

\begin{center}
{\huge \textbf{On the Determination of the Number of  \vskip.2cm Positive and Negative Polynomial Zeros  \vskip.35cm and Their Isolation}}

\vspace {10mm}
\noindent
{\large \bf Emil M. Prodanov} \vskip.5cm
{\it School of Mathematical Sciences, Technological University Dublin, Ireland,}
\vskip.1cm
{\it E-Mail: emil.prodanov@dit.ie} \\
\vskip1cm
\end{center}

\vskip4cm

\begin{abstract}
\noindent
A novel method with two variations is proposed with which the number of positive and negative zeros of a polynomial with real coefficients and degree $n$ can be restricted with significantly better determinacy than that provided by the Descartes rule of signs and also isolate quite successfully the zeros of the polynomial. The method relies on solving equations of degree smaller than that of the given polynomial. One can determine analytically the exact number of positive and negative zeros of a polynomial of degree up to and including five and also fully isolate the zeros of the polynomial analytically and with one of the variations of the method, one can analytically approach polynomials of degree up to and including nine by solving equations of degree no more than four. For polynomials of higher degree, either of the two variations of the method should be applied recursively. Full classification of the roots of the cubic equation, together with their isolation intervals, is presented. Numerous examples are given.
\end{abstract}

\newpage

\section{Introduction}

An algebraic equation with real co-efficients cannot have more positive real roots than the sequence of its co-efficients has variations of sign is the statement of Descartes' original rule of signs \cite{d} from 1637. Gauss showed \cite{g} in 1876 that the number of positive real roots (counted with their multiplicity) is, more precisely, either equal to the number of variations of signs in the sequence of the co-efficients, or is equal to the number of variations of signs in the sequence of the co-efficients reduced by an even number. \\
Many extensions of the Descartes rule have been proposed --- see \cite{marden}, where Marden has given a thorough summary of the results about polynomial roots. For example, if $I = (a, b)$ is an arbitrary open interval, the mapping $x \mapsto (ax+b)/(x+1)$ maps $(0, \infty$) bijectively onto $(a, b)$. Hence, if $\p$ is a polynomial of degree $n$, then the positive real zeros of $(1 + x)^n \, p_n[(ax+b)/(x+1)]$ correspond bijectively to the real zeros of $\p$ in $I$ \cite{m}. \\
Bounds on the zeros of polynomials were first presented by Lagrange \cite{l} and Cauchy \cite{c}. Ever since, the determination of the number of positive and negative roots of an equation, together with finding root bounds has been subject of intensive research. A more recent survey is provided by Pan \cite{pan}. Currently, the best root isolation techniques are subdivision methods with Descartes' rule terminating the recursion. \\
In this work, a novel method is proposed for the determination of the number of positive and negative zeros of a given polynomial $\p$ with real co-efficients and of degree $n$. The method also allows to find bounds on the zeros of $\p$. These bounds will not be on all of the zeros as a bulk, but, rather, if the bounds are not found individually, thus isolating each of the zeros, then not many of the zeros of the polynomial would be clumped into one isolation interval. All of this is achieved by considering the given polynomial $\p$ as a difference of two polynomials, the intersections of whose graphs gives the roots of $\p = 0$. The idea of the method is to extract information about the roots of the given polynomial by solving equations of degree lower than that of $\p$ --- in some sense by ``decomposing" $\p$ into its ingredients and studying the interaction between them. Different decompositions (further referred to as ``splits") yield different perspectives. Two splits are studied and illustrated in detail with numerous examples (with a different approach to one of the splits mentioned at the end). The first variation of the method splits the polynomial $\p$ by presenting it as a difference of two polynomials one of which is of degree $n$, but for which the origin is a zero of order $k < n$, while the other polynomial in the split is of degree $k-1$. For instance, polynomials of degree 9 can be split in the ``middle" with $k = 5$ in which case all  resulting equations that need to be solved are of degree 4 and their roots can be found analytically. This proves to be a very rich source of information about the zeros of this given polynomial of degree 9. This is illustrated with an example in which the negative zeros, together with one of the positive zeros of the example polynomial, have all been isolated, while the remaining two positive zeros are found to be within two bounds. All of the respective bounds are well within the bulk bounds, as found using the Lagrange and Cauchy formul\ae. The advantages over the Descartes rule of signs are also clearly demonstrated. If the degree of $\p$ is higher, then the method should be applied recursively. \\
The second variation of the method splits $\p$ by selecting the very first and the very last of its terms and ``propping" them against the remaining ones. This split is studied in minute detail and (almost) full classification of the roots of the cubic equation is presented, together with their isolation intervals and the criterion for the classification. The idea again relies on the ``interaction" between two, this time different, ``ingredients" of $\p$. One of these is a curve the graph of which passes through point $(0, 1)$, while the graph of the other one passes through the origin. By varying the only co-efficient of the former and by solving equations of degree $n-1$, one can easily find the values which would render the two graphs tangent to each other and also find the points at which this happens. Then simple comparison of the given co-efficient of the leading term to these values immediately determines the exact number of positive and negative roots and their isolation intervals for any equation of degree 5 or less. For polynomials of degree 6 or more, one again has to apply the method recursively. Such application is shown through an example with equation of degree 7. \\
To demonstrate how the proposed method works, it is considered on its own and no recourse is made to any of the known methods for determination of the number of positive and negative roots of polynomial equations or to any techniques for the isolation of their roots, except for comparative purposes only.

\section{The Method}

\n
Consider the equation
\b
\label{eq}
a_n \, x^n + a_{n-1} \, x^{n-1} + \ldots + a_1 \, x + a_0 = 0
\e
and write the corresponding polynomial $\p$ as
\b
\label{p1}
\p =  a_n \, x^n + a_{n-1} \, x^{n-1} + \ldots + a_1 \n x + a_0 = f_n(x) - g_{k-1}(x),
\e
where $0< k < n$ and:
\b
f_n(x) & = & x^k \, (a_n \, x^{n-k} + a_{n-1} \, x^{n-k-1}  + \ldots + a_{k+1} \, x + a_k) \equiv x^k \, F_{n-k}(x), \\
g_{k-1}(x) & = & - a_{k-1} \, x^{k-1} - a_{k-2} \, x^{k-2} - \ldots - a_1 \, x - a_0.
\e
The roots of the equation $p_n(x) = f_n(x) - g_{k-1}(x) = 0$ can be found as the abscissas of the intersection points of the graphs of the polynomials $f_n(x)$ and $g_{k-1}(x)$. The polynomial $f_n(x)$ has a zero $f_0$ of order $k$ at the origin and if $a_k \ne 0$ and $k$ is odd, it has a saddle there, while if $a_k \ne 0$ and $k$ is even, $f_n(x)$ has a minimum or a maximum at the origin. The remaining roots of $f_n(x) = 0$ are those of $F_{n-k}(x) = a_n x^{n-k} + a_{n-1} x^{n-k-1}  + \ldots + a_{k+1} x + a_k = 0$. The root $f_0 = 0$, together with the real roots $f_i$ and $g_i$ of the lower-degree equations $F_{n-k}(x) = 0$ and $g_{k-1}(x) = 0$, respectively, divide the abscissa into sub-intervals. The roots of the equation $\p = 0$ can exist only in those sub-intervals where $f_n(x)$ and $g_{k-1}(x)$ have the same signs and this also allows to count the number of positive and negative roots of the given equation $\p = 0$. When counting, one should keep in mind that the function $f_n(x)$ can have up to $n-1$ extremal points with the origin being an extremal point of order $k-1$, while the function $g_{k-1}(x)$ can have up to $k-2$ extremal points. This places an upper limit on the number of sign changes of the first derivatives of $f_n(x)$ and $g_{k-1}(x)$, that is, un upper limit on the ``turns" which the polynomials $f_n(x)$ and $g_{k-1}(x)$ can do, and this, in turn, puts an upper limit of the count of the possible intersection points between $f_n(x)$ and $g_{k-1}(x)$ in the various sub-intervals. If there is an odd number of roots of $f_n(x) = 0$ between two neighbouring roots of $g_{k-1}(x) = 0$ [or, vice versa, if there is an odd number of roots of $g_{k-1}(x) = 0$ between two neighbouring roots of $f_n(x) = 0$], then there is an odd number of roots of the equation $\p = 0$ between these two neighbouring roots of $g_{k-1}(x) = 0$ [or between the two neighbouring roots of $f_n(x) = 0$]. But if there is an even number of roots of $f_n(x) = 0$ between two neighbouring roots of $g_{k-1}(x) = 0$ (or vice versa), then the number of roots of $\p = 0$ between these two neighbouring roots of $g_{k-1}(x) = 0$ [or between the two neighbouring roots of $f_n(x) = 0$] is zero or some even number. In all cases, the end-points of these sub-intervals serve as root bounds and thus the number of positive and negative roots can be found with significantly higher determinacy than the Descartes rule of signs  provides. \\
All of the above is doable analytically for equations of degree up to and including 9 and is illustrated further with examples. In the case of $p_9(x)$, one will only have to solve two equations of degree four. \\
If one has a polynomial equation of degree 10 or higher, then the above procedure should be done recursively at the expense of reduced, but not at all exhausted, ability to determine root bounds and number of positive and negative roots. \\
To introduce a variation of the method, an assumption will be made: $\p$ is such that $0$ is not among the roots of the corresponding polynomial equation $\eq$, that is $a_0 \ne 0$. As the  determination of whether 0 is a root of an equation or not is absolutely straightforward, the case of a root being equal to zero will be of no interest for the analysis. In view of this, the co-efficient $a_0$ will be set equal to 1. It suffices to say that, should the equation $\p = 0$ has zero as root of order $m$, then the remaining non-zero roots of the equation $\p = 0$ can be found as the roots of the equation of degree $n-m$ given by $\tau_{n-m}(x) = p_n(x) / x^m = 0$. \\
One can consider an alternative split of the given polynomial $\p$ --- it can be written as the difference of two polynomials, each of which passes through a fixed point in the $(x, y)$--plane:
\b
\label{p2}
\p = a_n \, x^n + a_{n-1} \, x^{n-1} + \ldots + a_1 \, x + 1 =   \q - \r,
\e
where:
\b
\label{q}
\q & = & a_n \, x^n + 1, \\
\label{r}
\r & = & - a_{n-1} \, x^{n-1} - a_{n-2} \, x^{n-2} - \ldots  - a_1 \, x.
\e
The roots of the equation $p_n(x) = \q - \r = 0$ are found as the abscissas of the intersection points of the graphs of the polynomials $\q$ and $\r$. Regardless of its only co-efficient $a_n$, the polynomial $\q$ passes through point $(0, 1)$ --- the reason behind the choice of $a_0 = 1$, --- while the polynomial $\r$ passes through the origin (it has a zero there), regardless of the values of its co-efficients $a_{n-1}, a_{n-2}, \ldots , a_1$. \\
The method will be illustrated first for this split. \\
Considering the given $\p$, write $\alpha$ instead of the given coefficient $a_n$ and treat this $\alpha$ as undetermined. All other co-efficients $a_j, \,\, j = 1, 2, \ldots, n-1,$ are as they were given through the equation. Calculate next the discriminant $\Delta_n(\alpha)$ of the given polynomial $p_n(x)$. If this discriminant is zero, then the equation $\eq$ will have at least one repeated root. The equation $\Delta_n(\alpha) = 0$ is an equation in $\alpha$ of degree $n-1$. Denote the $n-1$ roots of the equation $\Delta_n(\alpha) = 0$ by $\alpha_1, \, \alpha_2 \, \ldots, \, \alpha_{n-1}$. Then, for the real roots of $\Delta_n(\alpha) = 0$, each of the equations
\b
\label{dv}
\alpha_j x^n + a_{n-1} x^{n-1} + \ldots + a_1 x + 1 = 0, \qquad \mbox{ $j \le n-1$},
\e
will have a root $\beta_j$ of order at least 2. If, in each of the above equations, $\alpha_j$ is perturbed slightly, so that $\Delta_n(\alpha_j)$ becomes negative for that perturbed $\alpha_j$, then the double real root $\beta_j$ will become a pair of complex conjugate roots. If, instead, the perturbation of $\alpha_j$ results in $\Delta_n(\alpha_j)$ becoming positive, then the real double root $\beta_j$ will bifurcate into two different real roots --- one on each side of $\beta_j$. \\
It should be noted that if $\alpha_j$ are all complex, then the equation $\eq$ cannot have a repeated root, namely, the equation $\eq$ has no real roots if it is of order $2m$ or has just one real root if it is of order $2m+1$. \\
Consider now $\q$ and $\r$. If the equation $p_n(x) = \q - \r = 0$ has a double root $\chi$, then the curves $\q$ and $\r$ will be tangent to each other at $\chi$, that is $q_n(\chi) = r_{n-1}(\chi)$, and, also, the tangents to the curves $\q$ and $\r$ will coincide at $\chi$, namely $q_n'(\chi) = r_{n-1}'(\chi)$. The latter allows to find
\b
a_1 = - 2 a_2 \chi - 3 a_3 \chi^2 - \ldots - n a_n \chi^{n-1}.
\e
Substituting into $q_n(\chi) = r_{n-1}(\chi)$ yields:
\b
- (n - 1) a_n \chi^{n} - (n - 2) a_{n - 1} \chi^{n-1} - \ldots - a_3 \chi^3  - a_2 \chi^2 + 1 = 0.
\e
Adding $(n-1)$ times $q_n(\chi) - r_{n-1}(\chi) = 0$ to the above results in the following equation:
\b
\label{chi}
a_{n-1} \chi^{n-1} + 2 a_{n-2} \chi^{n-2} + \ldots + (n-1) a_1 \chi + n = 0.
\e
The roots $\chi_j$ of this equation are the same as the double roots of the equations (\ref{dv}). Equation (\ref{chi}) is another equation of order one less than that of the original equation. \\
For an equation of degree up to and including 5, comparison of the given coefficient $a_5$ to the real numbers $\alpha_j$ from the obtained set $\,\, \alpha_1, \,\, \alpha_2, \,\, \alpha_3, \,\, \alpha_4\,\,$ allows not only to determine the exact number of positive and negative roots of the equation but to also isolate them. \\
For equation of degree 6 or higher, this variant of the method should be used recursively.

\section{Examples for the Split (\ref{p2})}

\n
The method will now be illustrated with examples (of increasing complexity) of polynomial equations of different orders. In the trivial case of $p_1(x) = 0$, that is $ax + 1 = 0$, the only root is determined by the intersection of the straight line $y = ax + 1$ with the abscissa $y = 0$. The root $x = - 1/a$ always exists. It is positive when $a < 0$ and negative when $a > 0$. The case of quadratic equation is also very simple --- see Figure 1 for the full classification.
\begin{figure}[!ht]
\centering
\subfloat[\scriptsize The case of $a > 0$ and $b \le 0$. When $a > 0$ and $b \ge 0$, the situation is analogical to the one shown here as there is axial symmetry (reflection) with respect to the ordinate (one replaces $b$ with $-b$). The points $x_{1,2}$ are the roots of the quadratic equation.  For fixed $a$, there is a value of $b$, say $\beta$, such that at point $x = \chi$, the graphs of $q_2(x)$ and $r_1(x)$ are tangent. The tangents to the graphs at point $x = \chi$ also coincide. Thus $q_2(x)\vert_{(x = \chi, b = \beta)} = r_1(x)\vert_{(x = \chi, b = \beta)}$ and $(d/dx) q_2(x)\vert_{(x = \chi, b = \beta)} = (d/dx) r_1(x)\vert_{(x = \chi, b = \beta)}$. That is, one has the two simultaneous equations $a\chi^2 + 1 = - \beta \chi$ and $2 a \chi = - \beta$, from which it is easily determined that $\beta = - 2 \sqrt{a}$ and $\chi = -\beta/2a = 1/\sqrt{a}$. The resulting root $\chi$ of the quadratic equation is double. This corresponds to vanishing discriminant $\Delta_2 = b^2 - 4a$. For $ - 2\sqrt{a} < b \le 0$ (i.e. when the discriminant $\Delta_2$ is negative), there are no real roots of the quadratic equation (this includes the depressed equation $ax^2 + 1 = 0$ for which $b = 0$). When $b < -2\sqrt{a}$, the double root bifurcates into two: $x_1$ and $x_2$ so that they fall on either side of the double root $\chi$. For $a > 0$ and $b < 0$, the roots, if they are real, are both positive. Thus, one of the roots of the quadratic equation $ax^2 + bx + 1 = 0$ is between the origin and $1/\sqrt{a}$. The other one is greater than $1/\sqrt{a}$. For $a > 0$ and $b < 0$ (the axially symmetric case), both roots, if they are real, are negative with one smaller than $-1/\sqrt{a}$ and the other --- between $-1/\sqrt{a}$ and the origin.  ]
{\label{F1a}\includegraphics[height=6cm, width=0.48\textwidth]{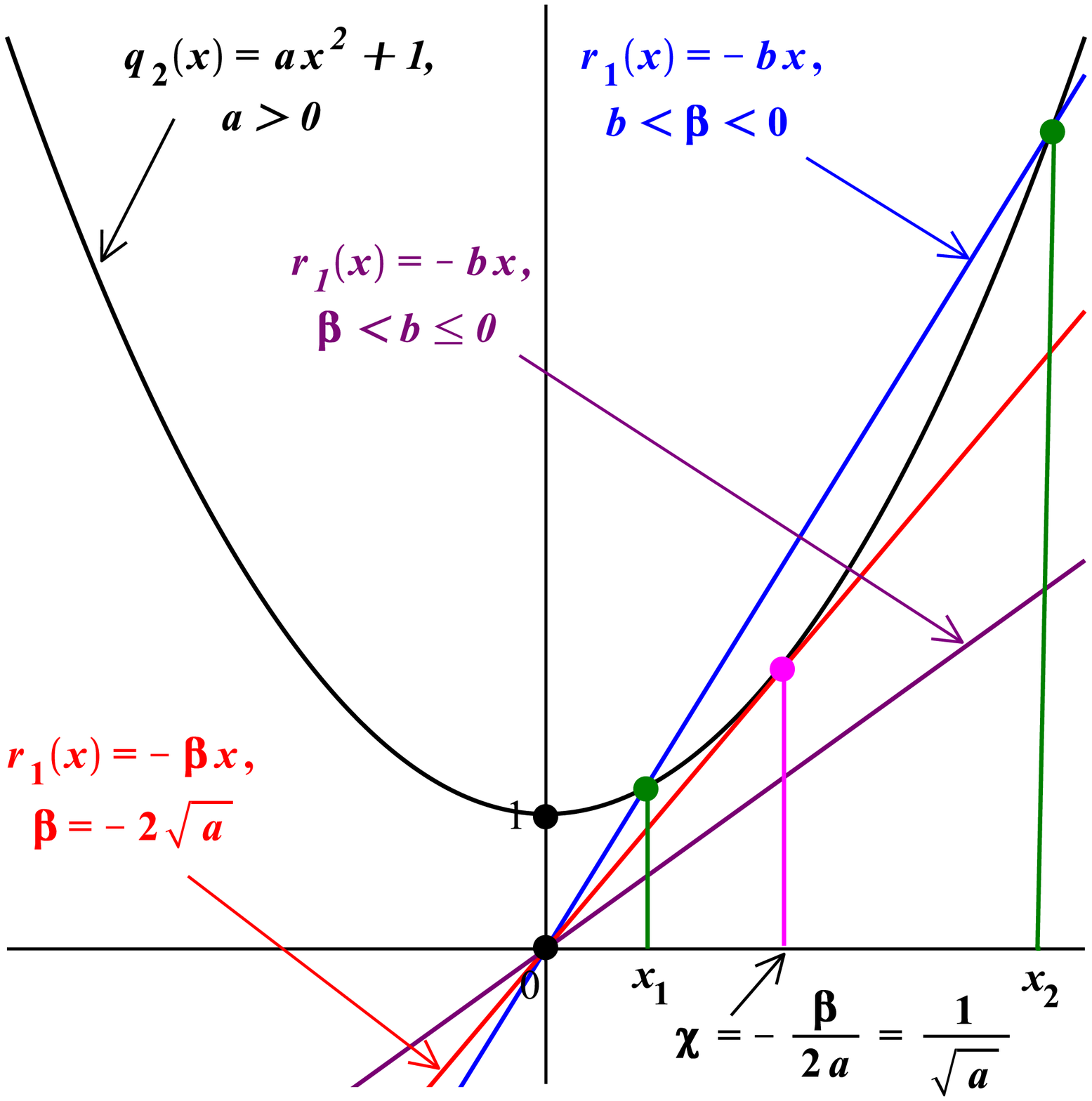}}
\quad
\subfloat[\scriptsize The case of $a < 0$ with: $b > 0$, $b = 0$ (corresponding to the depressed equation $ax^2+1=0$), and $b < 0$. The roots of the quadratic equation are real in each of the three subcases. One root is always positive, the other root is always negative. The roots of the suppressed equation are $\pm 1/\sqrt{-a}$ (recall, $a < 0$ now). Thus, the graph of the quadratic polynomial $q_2(x) = ax^2+1$ always goes through point $(0,1)$ and also through  points $(\pm 1/\sqrt{-a}, 0)$. When $b < 0$, the bigger root of the quadratic equation is between $0$ and $1/\sqrt{-a}$, while the smaller root is less than $-1/\sqrt{-a}$. When $b > 0$, the situation is symmetric: the smaller root is between $-1/\sqrt{-a}$ and $0$, while the bigger root is greater than $1/\sqrt{-a}$. With $a$ fixed and $b \to \infty$, the negative root tends to zero and the positive root tends to $\infty$.  With $a$ fixed and $b \to -\infty$, the negative root tends to $-\infty$ the positive root tends to $0$. Clearly, with $a$ fixed and $b \to 0^{\pm}$, the roots tend to those of the suppressed equation. With $b$ being a fixed positive number and $a \to 0^{-}$, the positive root tends to $+ \infty$ and the negative root tends to $-1/b$. With $b$ being a fixed negative number and $a \to 0^{-}$, the negative root tends to $- \infty$ and the positive root tends to $1/b$. When $b$ is a fixed number (positive, negative, or zero) and $a \to -\infty$, the roots tend to $0$ from either side, regardless of the sign of $b$. ]
{\label{F1b}\includegraphics[height=6cm,width=0.48\textwidth]{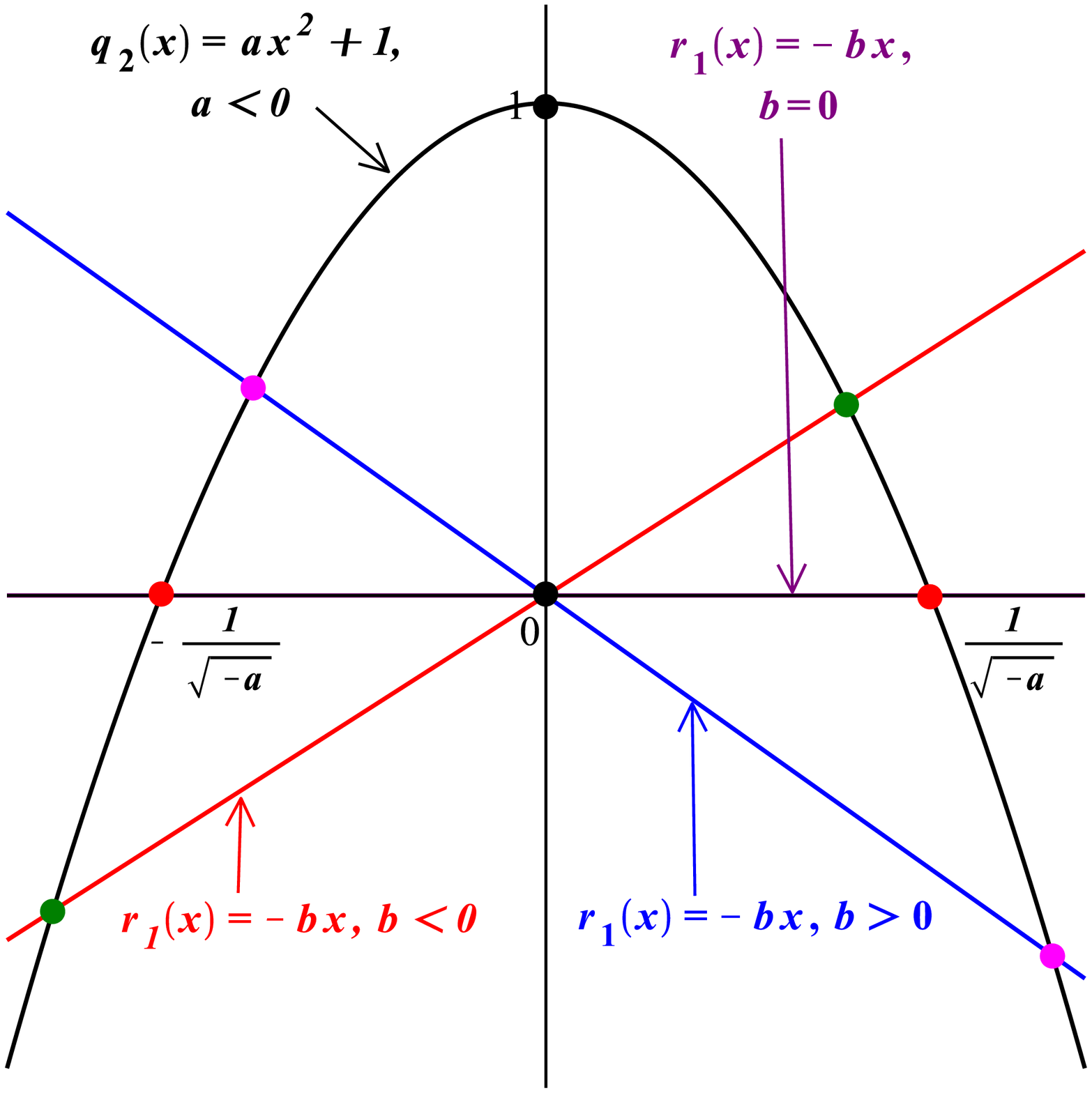}}
\caption{\footnotesize{The quadratic equation $ax^2 + 1 = - bx$.}}
\label{Figure1}
\end{figure}

\subsection{Cubic Equations}

\n
Consider a cubic equation which does not have a zero root. Without loss of generality, such equation can be written as $p_3(x) = ax^3 + bx^2 + cx + 1 = 0$ or as either one of the following two splits: $p_3(x) = q_3(x) - r_2(x) = 0$ with $q_3(x) = a x^3 + 1$ and $r_2(x) = - b x^2 - c x$, or $p_3(x) = \psi_3(x) - \varphi_1(x) = 0$ with $\psi_3(x) = a x^3 + b x^2$ and $\varphi_1(x) = - c x - 1$. That is:
\b
\label{q-r}
a x^3 + 1 & \!\! = \!\! & - b x^2 - c x, \\
\label{psi-phi}
a x^3 + b x^2 & \!\! = \!\! & - c x - 1.
\e
Provided that $b \ne 0$, the graph of the quadratic polynomial $r_2(x) = - b x^2 - c x$ passes through the origin and also through point $- c/b$ from the abscissa. The graph of the cubic polynomial $q_3(x) = a x^3 + 1$ passes through point $(0, 1)$ and also through point $-\sqrt[3]{1/a}$ from the abscissa. The cubic polynomial $\psi_3(x) = a x^3 + b x^2$ has a double root at zero. If $b$ is positive, the graph has a minimum at the origin, if $b$ is negative, the graph has a maximum at the origin. The other zero of the cubic polynomial $\psi_3(x) = a x^3 + b x^2$ is at $-b/a$. The graph of the polynomial $\varphi_1(x) = - c x - 1$ is a straight line passing through point $(0, -1)$ with slope $-c$. \\

\subsubsection{The Depressed Cubic Equation}

\n
Firstly, the situation of the depressed cubic equation, i.e. equation with $b = 0$, will be considered. In this case, the two splits become equivalent: the graphs of the pair in the split (\ref{q-r}) are the same as the graphs of the pair in the split (\ref{psi-phi}), but shifted vertically by one unit. To find for which $\alpha$ the equation $\alpha x^3 + c x + 1 = 0$ would have a double root, consider the discriminant $- \alpha (4c^3 + 27 \alpha)$. This vanishes when $\alpha$ is either $0$ or $-4c^3/27$.  The cubic equation $\alpha x^3 + cx + 1 = 0$ with $\alpha = -4c^3/27$ has root $\chi_0 = 3/c$ and a double root $\chi = -3/(2c)$. \\
One can then immediately determine the number of positive and negative roots of the original equation  $a x^3 + cx+1 = 0$  and also localise them as follows (see Figure 2).
\begin{figure}[!ht]
\centering
{\label{F2}\includegraphics[height=8cm, width=0.58\textwidth]{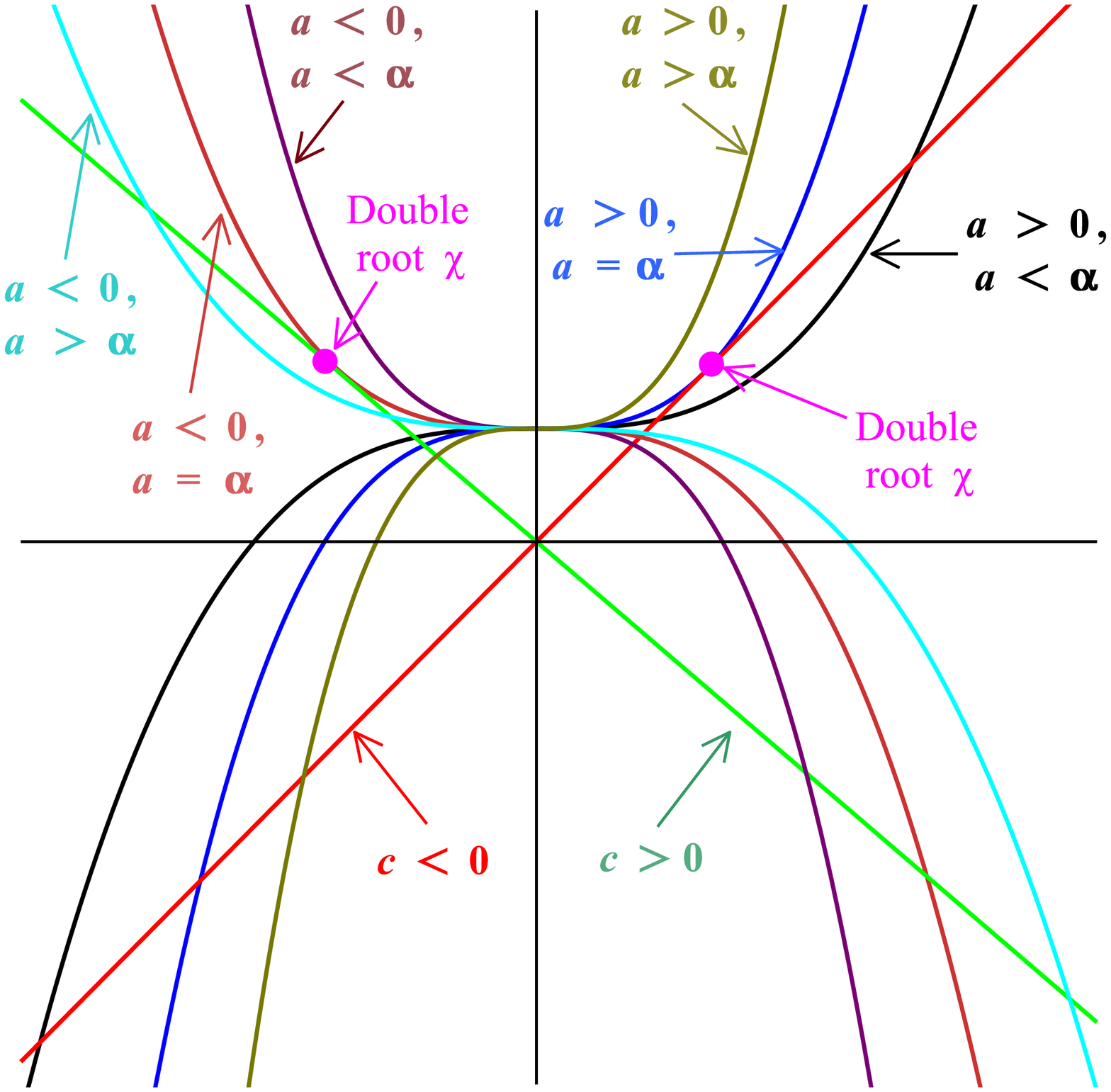}}
\caption{\footnotesize{The depressed cubic equation $ax^3 + 1 = -cx$. Equation $\alpha x^3 + 1 = -cx$ with $\alpha = -4c^3/27$ has root $\chi_0 = 3/c$ and a double root $\chi = -3/(2c)$. Comparing the given $a$ to $\alpha$ allows the immediate determination of the exact number of positive and negative roots of the equation and also their localisation.}}
\label{Figure2}
\end{figure}
\n
If $a$ and $c$ are both positive, then the equation has one negative root $x_0$ located between the intersection point of the graph of $ax^3+1$ with the abscissa and the origin. Namely, $\sqrt[3]{-1/a} < x_0 < 0$. The other two roots are complex. \\
If $a$ is positive and $c$ is negative, then the equation has one negative root $x_0$ located to the left of the intersection point of the graph of $ax^3+1$ with the abscissa. Namely, $x_0 < \sqrt[3]{-1/a}$. If, further, $a > \alpha$, then the equation has two complex roots. If $a = \alpha$, the equation has an additional double root given by $\chi = -3/2c > 0$. If, further, $a < \alpha$, then the equation has, additionally, two positive roots: one between 0 and $\chi = -3/2c$, the other --- greater than $\chi = -3/2c$. One only needs to compare $a$, given through the equation, to $\alpha = -4c^3/27$ for which the cubic discriminant vanishes. \\
If $a$ is negative and $c$ is positive, the situation is symmetric (with respect to the ordinate) to the situation of $a$ positive and $c$ negative. This time, the equation has one positive root $x_0$ located to the right of the intersection point of the graph of $ax^3+1$ with the abscissa. Namely, $x_0 > \sqrt[3]{-1/a}$. If, further, $a < \alpha$, then the equation has no other roots. If $a = \alpha$, the equation has an additional double root given by $\chi = -3/2c < 0$. Finally, if $a > \alpha$, then the equation has, additionally, two negative roots: one between 0 and $\chi = -3/2c$, the other --- smaller than $\chi = -3/2c$. Again, one only needs to compare $a$, given through the equation, to $\alpha = -4c^3/27$ for which the cubic discriminant vanishes. \\
Finally, if both $a$ and $c$ are negative, then the equation has one positive root $x_0$ located between the origin and the intersection point of the graph of $ax^3+1$ with the abscissa, that is, $0 < x_0 < \sqrt[3]{-1/a} < x_0$. \\
As a numerical example of the above, consider the equation
\b
\label{orig}
3 x^3 - 4x + 1 = 0.
\e
The auxiliary equation
\b
\label{aux}
\alpha x^3 - 4x + 1 = 0
\e
has discriminant $ - \alpha (27 \alpha - 256)$. Clearly, when $\alpha = 256/27$, equation (\ref{aux}) will have a double root $\chi$. To find $\chi$, one needs to solve the equation (\ref{chi}):
\b
2 c \chi + 3 = 0,
\e
with $c = -4$. Thus, $\chi = 3/8 = 0.375$. \\
As the given value of $a$ is 3 and as $a = 3 < \alpha = 256/27 \approx 9.481,$ the given equation (\ref{orig}) has two positive roots: $x_2$ which is between 0 and $\chi = 0.375$, and $x_3$ which is greater than $\chi = 0.375$. The equation also has a negative root $x_1$ to the left of the intersection point of the graph of $3 x^3 + 1$ and the abscissa, that is $x_1 < - \sqrt[3]{1/3} \approx -0.693$. \\
The roots of equation (\ref{orig}) are: $x_1 = -1.264,\,\, x_2 = 0.264, \,\, x_3 = 1.$

\subsubsection{Full Cubic Equation}

\n
If one considers next the cubic equation $ax^3 + bx^2 + 1 = 0$, the situation will not turn out to be qualitatively different from the one of the ``full" equation $ax^3 + bx^2 + cx + 1 = 0$, where all co-efficients are different from zero, and the latter is the equation to be considered next. \\
The discriminant of the cubic equation
\b
\label{cubic}
a x^3 + b x^2 + c x + 1 = 0
\e
is given by
\b
\label{d3}
\Delta_3 = - 27 a^2 + 2c (9b - 2c^2) a + b^2 (c^2 - 4b).
\e
Setting $\Delta_3 = 0$ and interpreting $b$ and $c$ as parameters, one gets a quadratic equation for $a$ with roots
\b
\label{alfa12}
\alpha_{1,2} = - \frac{2}{27} \, c^3 + \frac{1}{3} \, bc \pm \frac{2}{27} (c^2 - 3b)^{\frac{3}{2}}.
\e
These are real, i.e. the discriminant $\Delta_3$ can be zero, only when $c^2 > 3b$. Let $\alpha_1$ denote the bigger root. \\
If, further, $c^2 > 4b$, then the free term in the quadratic equation $\Delta_3 = 0$ will be positive and, according to Vi\`ete's formula, $\alpha_1$ and $\alpha_2$ will have different signs.  \\
The two border cases are $c^2 = 4b$, in which case the roots (\ref{alfa12}) are $0$ and $c^3/54$, and $c^2 = 3b$, in which case equation $\Delta_3 = 0$ has a double root $c^3/27$. \\
For the case of a general equation of degree 3, equation (\ref{chi}) becomes:
\b
b x^2 + 2 c x + 3 = 0.
\e
The roots of this equation are
\b
\chi_{1,2} = \frac{1}{b} \left( -c \pm \sqrt{c^2 - 3b} \right).
\e
At point $x = \chi_{1,2}$, the curve $\alpha_{1,2} x^3 + 1$ is tangent to the curve $- b x^2 - c x$ and the tangents to the graphs to each of the  curves coincide at that point. \\
To determine the exact number of positive and negative roots and to also localise the roots, one has to compare the given $a$ with the values of $\alpha_1$ and $\alpha_2$. \\
Depending on the signs of the three parameters $a$, $b$, and $c$, there are eight cases to be analysed. Four will be considered in detail, the analysis for each of the remaining four cases can be easily inferred afterwards. \\ \\
\subsubsubsection{}
\vskip-0.73cm
\paragraph{\hskip1.5cm The Case of $\bm{a > 0, \,\, b > 0$}, and $\bm{c >0}$ \\ \\}
$\hskip-.53cm$ This is the most complicated case. There are three sub-cases.
\subparagraph{\hskip-.6cm 2.1.2.1.1 The Sub-case of $\bm{\,\,\, -c/b < - \sqrt[3]{\frac{1}{a}}}$ \\ \\}
$\!\!\!\!\!\!\!\!$ In this sub-case, the negative root of $-bx^2 - cx$ is to the left of the point where $a x^3 + 1$ intersects the abscissa --- see Figure 3a. \\
If one further has $c^2 \ge 4b$, then $\alpha_1 > 0$ and $\alpha_2 \le 0$. The non-positive root $\alpha_2$ is not relevant to the analysis as $a$, given through the equation, is positive in this case (the curve $\alpha_2 x^3 + 1$ with $\alpha_2 \le 0$ is also tangent to the curve $-bx^2 - cx$, but the given equation $ax^3 + 1 = - bx^2 - cx$ has $a > 0$). Equation (\ref{chi}) has root $\chi_1 = (c/b) \left( - 1 - \sqrt{1 - 3b/c^2} \right) < -c/b < 0$ and root $\chi_2$, which is associated with $\alpha_2$. Thus, the double root $\chi$ of equation $\alpha_2 x^3 + 1 = - b x^2 - c x$ will be given by $\chi_1$, while the other root (using Vi\`ete's formula) will be $\chi_0 = - 2 \chi - b/\alpha_1$. \\
When $c^2 = 4b$, the curve $a x^3 + 1$ with $a = 0$ is tangent to $-bx^2 - cx$ at the maximum $c^2/(4b)$ of $-bx^2 - cx$ which occurs at $-c/(2b)$. When $c^2 > 4b$, the curve $a x^3 + 1$ with $a = 0$ will intersect the curve $-bx^2 - cx$ between $-c/(2b)$ and the origin at point $\sigma$ which is the bigger root of the quadratic equation $- bx^2 - cx = -1$, that is $ -c/(2b) < \sigma = [-c/(2b)] (1 - \sqrt{1-4b/c^2}) < 0$. Thus, for any $a$ greater than 0, provided that $-c/b < - \sqrt[3]{1/a}$, the curve $ax^3 + 1$ will intersect the curve $- bx^2 - cx$, for which $c^2 \ge 4b$, once between $\sigma$ and 0. \\
As a sub-case with $-c/b < - \sqrt[3]{1/a}$ is being studied, one has $a > (b/c)^3$. If $(b/c)^3 < a < \alpha_1$, the cubic equation will have two negative roots $x_1$ and $x_2$, such that $x_1 < \chi$ and $\chi < x_2 < -c/b$, where $\chi = (c/b) \left( - 1 - \sqrt{1 - 3b/c^2} \right)$ and $-c/b$ is the smaller root of $-bx^2 - cx =0$, and another negative root $x_3$ between $\sigma$ and the origin, i.e. $\sigma < x_3 < 0$. If $a \to \alpha_1$ from below, then the roots $x_1$ and $x_2$ will tend to $\chi$ from either side until they coalesce at the double root $\chi$ when $a = \alpha_1$. When $a > \alpha_1$, there will be a negative root between $\sigma$ and $0$ and two complex roots. \\
Consider as numerical examples for equation with $a > 0$, $b > 0$, $c > 0$, $c^2 \ge 4b$, and $-c/b < - \sqrt[3]{1/a}$ the following two equations.
\b
x^3 + 2x^2 + 3x + 1 = 0.
\e
The roots $\alpha_{1,2}$ are $\pm 2 \sqrt{3} / 9$. The relevant one is $\alpha_1 = 2 \sqrt{3} / 9 \approx 0.385$. The roots $\chi_{1,2}$ are $-(3/2) \pm \sqrt{3}/2$. The one of interest is $\chi = -(3/2) \pm \sqrt{3}/2 \approx - 2.366$. For this equation one has $\sigma = [-c/(2b)] (1 - \sqrt{1-4b/c^2}) = -1/2$. As $a = 1 > \alpha_1 \approx 0.385$, the given equation has one negative root between $\sigma = - 0.5$ and 0 and two complex roots. Indeed, the roots are: $-0.430$ and $-0.785 \pm  1.307 \, i$ approximately. \\
For the equation
\b
\frac{1}{3} \, x^3 + 2x^2 + 3x + 1 = 0
\e
one again has $\alpha_1 = 2 \sqrt{3} / 9 \approx 0.385, \,\,$ $\chi = -(3/2) \pm \sqrt{3}/2 \approx - 2.366,$ and $\sigma = -1/2$. This time $a = 1/3 < \alpha_1 \approx 0.385$. Therefore the given equation has three negative roots: one smaller than $\chi \approx -2.366$, another one between $\chi \approx -2.366$ and $-c/b = -1.5$, and the third one between $\sigma = -0.5$ and the origin. Indeed, the roots are: $-3.879, \,\,  -1.653,$ and  $-0.468$ approximately. \\
If one has $3b \le c^2 < 4b$ instead of $c^2 \ge 4b$, the roots $\alpha_{1,2}$ are both positive (with $\alpha_1 = \alpha_2$ when $c^2 = 3b$, otherwise  $\alpha_1 > \alpha_2$). The curve $\alpha_{1} x^3 + 1$ will be tangent to the curve $-bx^2 - cx$ at point $\chi_1 \le -c/b$, while curve $\alpha_2 x^3 + 1$ will be tangent to the curve $-bx^2 - cx$ at point $-(c/b) \le \chi_2 < -c/(2b)$. The points $\chi_{1,2}$ are the roots of the equation $bx^2 + 2cx+3=0$, that is $\chi_{1,2} = (c/b) \left( - 1 \pm \sqrt{1 - 3b/c^2} \right)$ --- see Figure 3b. If $c^2 = 3b$,
then $\alpha_1$ and $\alpha_2$ coalesce to $c^3/27$ while $\chi_1$ and $\chi_2$ coalesce to $-c/b = -3/c$. Then, the roots of the cubic equation $ax^3 + bx^2 + cx + 1 = 0$ will be as follows. If the given $a$ is such that $a > \alpha_1$, the cubic equation will have a negative root between $-\sqrt[3]{1/a}$ and 0 and two complex roots. If $a = \alpha_{1}$, there will be a negative root between $-\sqrt[3]{1/\alpha_1}$ and 0 and a negative double root at $\chi_{1}$. If $\alpha_2 < a < \alpha_1$, the cubic equation will have three negative roots $x_{1,2,3}$ such that $x_1 < \chi_1, \,\,\, \chi_1 < x_2 < \chi_2, \,$ and $ \chi_2 < x_3 < 0$. If $a = \alpha_{2}$, there will be a negative root smaller than $-\sqrt[3]{1/\alpha_2}$ and a negative double root at $\chi_{2}$. Finally, if $a < \alpha_2$ the cubic equation will have a negative root smaller than $-\sqrt[3]{1/a}$ and two complex roots. \\
As the current sub-case has the restriction $-c/b < - \sqrt[3]{1/a}$, one can only have $a > (b/c)^3$. Thus, when $(b/c)^3 < a  < \alpha_1$, the cubic equation will have two negative roots $x_1$ and $x_2$, such that $x_1 < \chi_1$ and $\chi_1 < x_2 < -c/b$, and another negative root $x_3$ between the intersection point of $ax^3 + 1$ with the abscissa and the origin: $-\sqrt[3]{1/a} < x_3 < 0$. When $a > \alpha_1$, the cubic equation will have one negative root between $-\sqrt[3]{1/a}$ and the origin together with two complex roots. \\
As numerical examples for equation with $a > 0$, $b > 0$, $c > 0$, $3b \le c^2 < 4b$, and $-c/b < - \sqrt[3]{1/a}$ consider the following two equations.
\b
x^3 + \frac{5}{3}x^2 + \frac{23}{10}x + 1 = 0.
\e
The roots $\alpha_{1,2}$ are $(5083 \pm 29 \sqrt{29})/13500$. That is, $\alpha_1 \approx 0.388$ and $\alpha_2 \approx 0.365$. The roots $\chi_{1,2}$ are $(-69 \pm 3 \sqrt{29})/50$, namely $\chi_1 \approx -1.703$ and $\chi_2 \approx -1.057$. As in the given equation $a = 1$, which is greater than the bigger root $\alpha_1 \approx 0.388$, the equation has a negative root between $-\sqrt[3]{1/a}$ and the origin, i.e. between $-1$ and $0$, and
two complex roots. Indeed, the roots of the equation are:  $-0.603$ and  $-0.532 \pm 1.173 \, i$ approximately. \\
Next, for the equation
\b
\frac{385}{1000} x^3 + \frac{5}{3}x^2 + \frac{23}{10}x + 1 = 0
\e
the roots $\alpha_{1,2}$ are the same: $\alpha_1 \approx 0.388$ and $\alpha_2 \approx 0.365$. The roots $\chi_{1,2}$ are also the same: $\chi_1 \approx -1.703$ and $\chi_2 \approx -1.057$. The smaller root of $-bx^2 - cx = 0$ is $-c/b = -1.38$ --- between $\chi_1$ and $\chi_2$, as expected. Also in this case: $(b/c)^3 \approx 0.381$. As the given $a$ is $0.385$, one has $\alpha_2 < (b/c)^3 < a  < \alpha_1$. Therefore, the roots of this cubic equation must be all negative and such that: one is smaller than $\chi_1 \approx -1.703$; another one is between $\chi_1 \approx -1.703$ and $- c/b = -1.38$; and the thrid one is between $-\sqrt[3]{1/a} \approx -1.375$ and 0. Indeed, the roots are $- 1.938, \,\, -1.494,$ and $ -0.897$ approximately. \\
Finally, when $c^2 < 3b$ in the sub-case of $-c/b < - \sqrt[3]{1/a}$, then the cubic equation will have one negative root between $- \sqrt[3]{1/a}$ and 0 and two complex roots. This is illustrated by
\b
2 x^3 + x^2 + x + 1 = 0.
\e
The model predicts a negative root between $- \sqrt[3]{1/a} \approx -0.794$ and 0 and two complex roots. Indeed, one has $x_1 = -0.739$ and $x_{2,3} =  0.119 \pm 0.814\,i$.

\begin{figure}[!ht]
\centering
\subfloat[\scriptsize For the sub-case of $-c/b < (-a)^{-1/3}$, the smaller root $-c/b$ of $-bx^2-cx$ is to the left of the intersection of $ax^3 + 1$ with the abscissa. If, further, $c^2 \ge 4 b$, the quadratic equation $\Delta_3 = 0$ will have roots $\alpha_1 > 0$ and $\alpha_2 \le 0$ (the latter bears no relevance as $a$ is considered positive in this case). Equation $\alpha_1 x^3 + 1 = - bx^2 - cx$ has a double negative root $\chi = (c/b) \left( - 1 - \sqrt{1 - 3b/c^2} \right) < -c/b $ and another negative root $\chi_0 = -2 \chi - b/\alpha_1$. Then, given the restriction $-c/b < - a^{-1/3}$ of this sub-case, i.e. $a > (b/c)^3$, one can have either: $(b/c)^3 < a < \alpha_1$ which leads to the cubic equation having two negative roots $x_1$ and $x_2$, such that $x_1 < \chi$ and $\chi < x_2 < -c/b$, together with another negative root $x_3$ between the origin and the bigger root $\sigma = \left(-c/(2b)\right) (1 - \sqrt{1-4b/c^2}) < 0$ of the equation $-bx^2 - cx = 1$ (the intersection point of $ax^3 + 1$ with $a = 0$ and $-bx^2 - cx$ closer to the origin), i.e. $\sigma < x_3 < 0$; or one can have $a = \alpha_1$ in which case the cubic equation will have a negative root between $\sigma$ and $0$ and a double root at $\chi$; or one can have $a > \alpha_1$, in which case the cubic equation will have a negative root between $\sigma$ and $0$ and two complex roots. If, instead of $c^2 \ge 4 b$, one has $c^2 <  3b$ (see Figure 3b for the case of $3b <c^2 < 4b$), then the cubic discriminant $\Delta_3$ will be negative and the cubic equation will have a negative root between $-a^{-1/3}$ and $0$ and two complex roots. If $c^2 = 3b$, then the double root $\chi$ is at $-c/b = -3/c$ and $\alpha_1 = \alpha_2 = c^3/27$. ]
{\label{F3a}\includegraphics[height=6cm, width=0.48\textwidth]{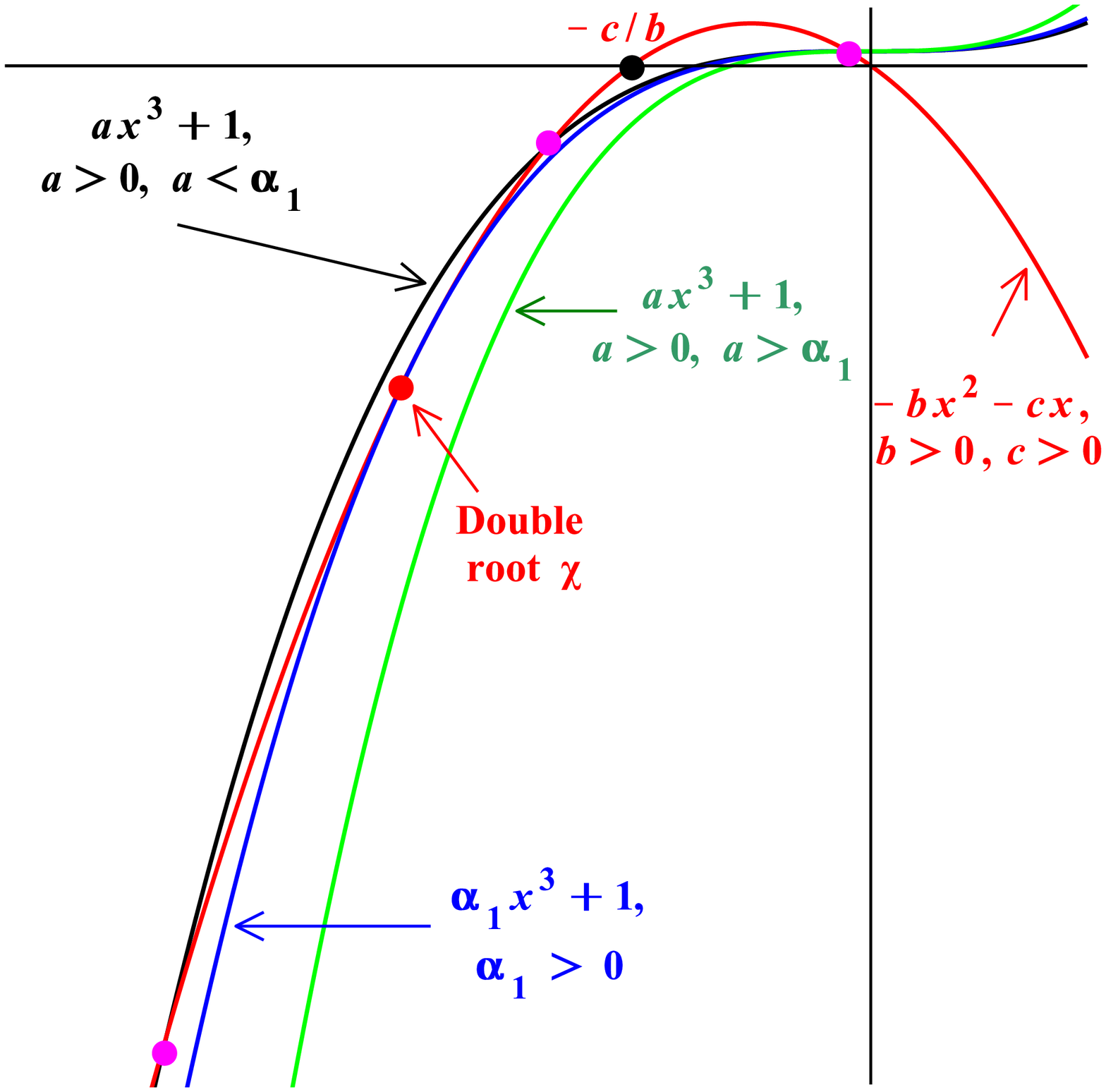}}
\quad
\subfloat[\scriptsize  If $3b < c^2 < 4b$ and with no further restrictions except $a > 0, \,\, b > 0$, and $c >0$, the roots $\alpha_{1,2}$ of the quadratic equation $\Delta_3 = 0$ will have the same signs. In this case the co-efficient $2c(9b - 2c^2)$ in the term linear in $a$ will be positive and thus both roots $\alpha_1$ and $\alpha_2$ will be positive (with $\alpha_1 > \alpha_2$). The curve $\alpha_{1} x^3 + 1$ will be tangent to the curve $-bx^2 - cx$ at point $\chi_1 < -c/b$, while curve $\alpha_2 x^3 + 1$ will be tangent to the curve $-bx^2 - cx$ at point $-c/b < \chi_2 < -c/(2b)$. The points $\chi_{1,2}$ are the roots of the equation $bx^2 + 2cx+3=0$, that is $\chi_{1,2} = (c/b) \left( - 1 \pm \sqrt{1 - 3b/c^2} \right)$. Then, the roots of the cubic equation $ax^3 + bx^2 + cx + 1 = 0$ will be as follows. If the given $a$ is such that $a > \alpha_1$, the cubic equation will have a negative root between $- a^{-1/3}$ and 0 and two complex roots. If $a = \alpha_{1}$, there will be a negative root between $- \alpha_1^{-1/3}$ and 0 and a negative double root at $\chi_{1}$. If $a$ is such that $\alpha_2 < a < \alpha_1$, the cubic equation will have three negative roots: $x_1 < \chi_1, \,\,\, \chi_1 < x_2 < \chi_2, \,$ and $ \chi_2 < x_3 < 0$. If $a = \alpha_{2}$, there will be a negative root smaller than $- \alpha_2^{-1/3}$ and a negative double root at $\chi_{2}$. Finally, if $a < \alpha_2$ the cubic equation will have a negative root smaller than $-a^{-1/3}$ and two complex roots. If $c^2 = 3b$, then $\alpha_1$ and $\alpha_2$ coalesce to $c^3/27$ while $\chi_1$ and $\chi_2$ coalesce to $-c/b = -3/c$. If $c^2 < 3b$, then $\chi_{1,2}$ are not among the reals. If $c^2 = 4b$, then $\alpha_2$ becomes 0 and $\chi_2 = -c/(2b)$ --- the point of the maximum of $- bx^2 - cx$. When $c^2 > 4b$, one has $\chi_2 > -c/(2b)$ and $\alpha_2 < 0$. ]
{\label{F3b}\includegraphics[height=6cm,width=0.48\textwidth]{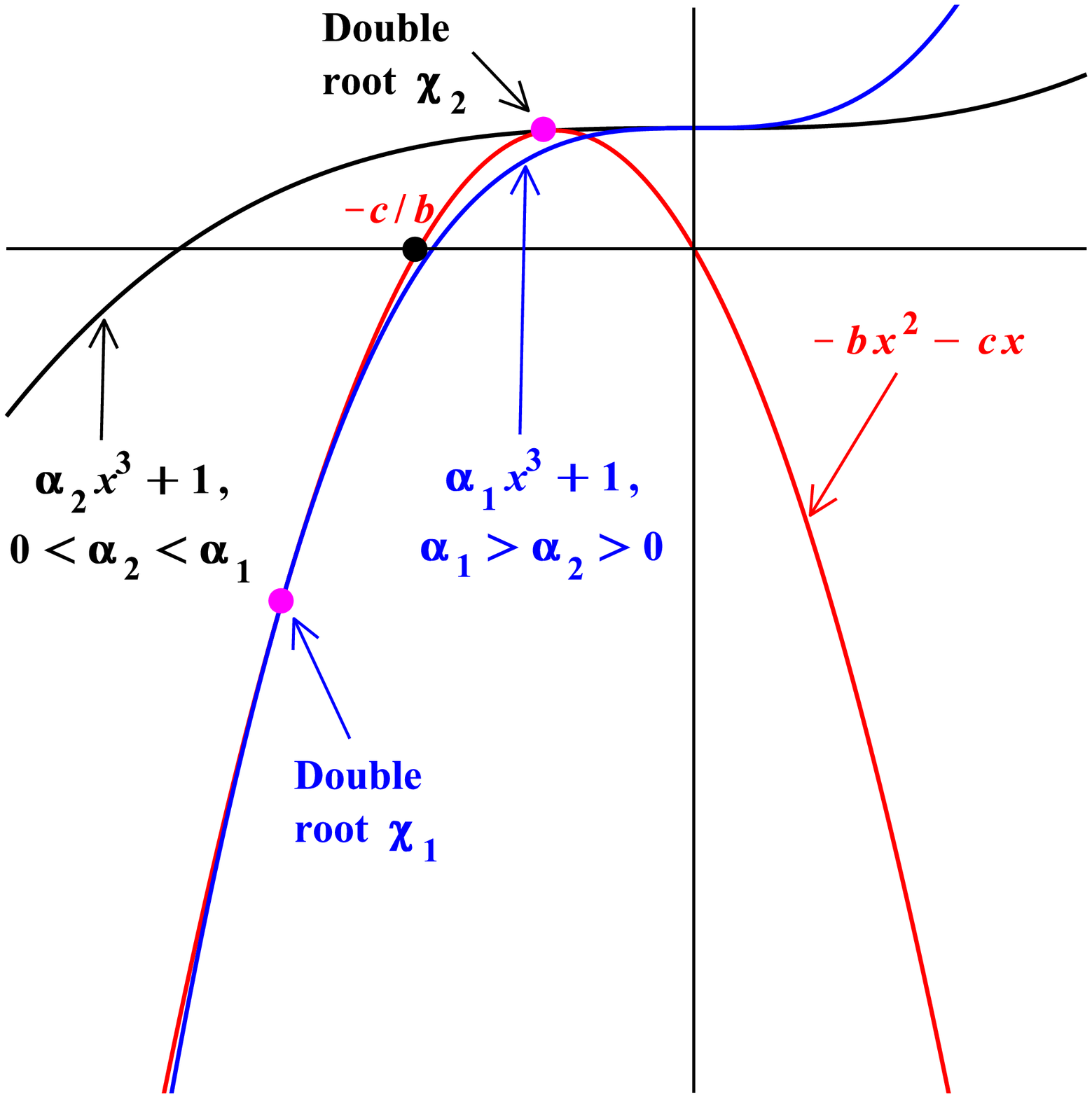}}
\caption{\footnotesize{The cubic equation $ax^3 + 1 = - b x^2 - c x$ with $a > 0$, $b > 0$, and $c > 0$.}}
\label{Figure3}
\end{figure}
\subparagraph{\hskip-.6cm 2.1.2.1.2 The Sub-case of $\bm{\,\,\, -c/b > - \sqrt[3]{\frac{1}{a}}}$ \\ \\}
$\!\!\!\!\!\!\!\!$ In this sub-case, the negative root of $-bx^2 - cx$ is between the point where $a x^3 + 1$ intersects the abscissa and the origin. \\
$\!\!\!\!\!\!\!\!\!$ If in this case one also has $c^2 > 4b$, then the curve $a x^3 + 1$ with $a = 0$ will intersect the curve $-bx^2 - cx$ between $-c/b$ and the origin at two points: the roots $\sigma_{1,2}$ of the quadratic equation $-bx^2 - cx = 1$, i.e. $\sigma_{1,2} = [-c/(2b)] (1 \pm \sqrt{1-4b/c^2}) < 0$ (with $\sigma_2 < \sigma_1$). These are on either side of $-c/(2b)$. Then for any $a$ greater than 0, given that $-c/b > - \sqrt[3]{1/a}$, the curve $ax^3 + 1$ will intersect the curve $- bx^2 - cx$ with $c^2 \ge 4b$, once between $-c/b$ and $\sigma_2$ and one more time  between $\sigma_1$ and the origin. There will be a third intersection to the left of $- \sqrt[3]{1/a}$. Therefore, the cubic equation will have 3 negative roots, the biggest of which will be between $\sigma_1$ and the origin, the middle one will be between $-c/b$ and $\sigma_2$, and the smallest one will be to the left of $-\sqrt[3]{1/a}$. A double root $\chi$ cannot exist when $c^2 \ge 4b$. \\
Equation with $a > 0$, $b > 0$, $c > 0$, $c^2 \ge 4b$, and $-c/b > - \sqrt[3]{1/a}$ can be illustrated with the following numerical example:
\b
\frac{1}{30} x^3 + x^2 + 3 x + 1 = 0.
\e
The roots are: $x_1 \approx -26.667, \,\,\, x_2 \approx -2.952$, and $x_3 \approx -0.381$ and within their predicted bounds: the biggest one is between $\sigma_1 \approx -0.382$ and the origin, the middle one is between $- c/b = -3$ and $\sigma_2 \approx -2.618$, and the third one is less than $-\sqrt[3]{1/a} \approx -3.107$. \\
Next, if one has $3b \le c^2 < 4b$ instead of $c^2 \ge 4b$, the situation on Figure 3c applies. In view of the restriction $-c/b > - \sqrt[3]{1/a}$, one can have: either $a < \alpha_2$, or $a = \alpha_2$, or $\alpha_2 < a < (b/c)^3$. In the first case, the cubic equation will have a negative root smaller than $- \sqrt[3]{1/a}$ and two complex roots. In the second case, the cubic equation will have a negative root smaller than $- \sqrt[3]{1/a}$ and a double negative root ar $\chi_2$. In the third case, the cubic equation will have three negative roots: $x_1 < \chi_1, \,\,\, \chi_1 < x_2 < \chi_2, \,$ and $ \chi_2 < x_3 < 0$. (As in the previous sub-case, if $c^2 = 3b$, then $\alpha_1$ and $\alpha_2$ coalesce to $c^3/27$ while $\chi_1$ and $\chi_2$ coalesce to $-c/b = -3/c$.) \\
Three numerical examples for equation with $a > 0$, $b > 0$, $c > 0$, $3b \le c^2 < 4b$, and $-c/b > - \sqrt[3]{1/a}$ are given. The first one is:
\b
\frac{2}{5} x^3 + 2x^2 + \frac{5}{2}x + 1 = 0.
\e
For this equation, the roots $\alpha_{1,2}$ are $14/27$ and $1/2$, i.e. $\alpha_1 \approx 0.519,\, \, \alpha_2 = 0.5$. The respective loci of the double roots are $\chi_1 = -3/2$ and $\chi_2 = - 1$. In this case, the given $a = 2/5 = 0.4$ is a number smaller than $\alpha_2 = 0.5$. Therefore, the equation has a negative root smaller than  $- \sqrt[3]{1/a} \approx - 1.357$ and two complex roots. The roots of the equation are $x_1 = -3.362$, and $x_{2,3} = -0.819 \pm 0.270\,i$. \\
Next, for the equation
\b
\frac{1}{2} x^3 + 2x^2 + \frac{5}{2}x + 1 = 0
\e
one has $\alpha_1 \approx 0.519$ and $\alpha_2 = 0.5$ with $\chi_1 = -3/2$ and $\chi_2 = - 1$. The given $a = 1/2$ is exactly equal to $\alpha_2 = 0.5$. Thus, there must be a negative root smaller than $- \sqrt[3]{1/a} \approx -1.256$ and a double root at $\chi_2 = -1$. The roots of this equation are $-2,\,\, -1,$ and $-1$ --- exactly in the predicted bounds. \\
In the third example, equation
\b
\frac{51}{100} x^3 + 2x^2 + \frac{5}{2}x + 1 = 0
\e
also has $\alpha_1 \approx 0.519$ and $\alpha_2 = 0.5$ with $\chi_1 = -3/2$ and $\chi_2 = - 1$. In this case, the given $a = 51/100$ is a number between $\alpha_2 = 0.5$ and $(b/c)^3 = 0.512$. The model predicts that the equation will have three negative roots such that the smallest one is smaller than $\chi_1 = -1.5$, the middle one is between $\chi_1 = -1.5$ and $\chi_2 = -1$, and the third one is between $\chi_2 = -1$ and 0. Indeed, the roots are: $-1.821,\,\, -1.213,$ and $-0.888$ approximately. \\
Finally, when $c^2 < 3b$ in the sub-case of $-c/b > - \sqrt[3]{1/a}$, then there is one negative root smaller than $- \sqrt[3]{1/a}$ and two complex roots. This is illustrated by the equation
\b
\frac{1}{2} x^3 + x^2 + x + 1 = 0.
\e
The negative root is $x_1 \approx -1.544 < - \sqrt[3]{1/a} \approx - 1.260$ and the other two roots are complex: $x_{2,3} \approx -0.228 \pm 0.115\, i$ --- as predicted.

 \subparagraph{\hskip-.6cm 2.1.2.1.3 The Sub-case of $\bm{\,\,\, -c/b = - \sqrt[3]{\frac{1}{a}}}$ \\ \\}
$\!\!\!\!\!\!\!\!\!\!$ When $a = b^3/c^3$, the cubic discriminant is $\Delta_3 = (9b^2/c^6) (c^2+b)(c^2-3b)$. It is negative when $c^2 < 3b$ and in this case the cubic equation will have a negative root given by $x = -c/b = \sqrt[3]{-1/a}$ together with two complex roots. If the discriminant $\Delta_3$ is positive (i.e. $c^2 > 3b$), the cubic equation will have three negative roots given by: $x_1 < -c/b = \sqrt[3]{-1/a}$, $x_2 = -c/b = \sqrt[3]{-1/a}$, and $x_3 > -c/b = \sqrt[3]{-1/a}$. Finally, when the discriminant $\Delta_3$ is zero, one will have $b = c^2/3$. Given that $a = b^3/c^3$, one gets that $a = c^3/27$. Thus, the cubic equation $c^3 x^3/27 + c^2 x^2/3 + cx + 1 = 0$ will have a triple negative root at $-3/c$. \\
The following three examples illustrate all possibilities for this sub-case. \\
Firstly, the equation
\b
x^3 + x^2 + x + 1 = 0,
\e
which is in the category $c^2 < 3b$ with $a = (b/c)^3$, should have a negative root given by $x = -c/b = -1$ together with two complex roots.
This is the case indeed, as the roots are: $x_1 = 1, \,\, x_{2,3} = \pm i$. \\
Next, equation
\b
\frac{8}{27} x^3 + 2x^2 + 3x + 1 = 0
\e
is with $c^2 > 3b$ and, since $a = (b/c)^3$, the roots are all negative: one is smaller than $-c/b = -1.5$, another one is equal to $-c/b = -1.5$ and the third one is between $-c/b = -1.5$ and 0. The roots of the equation are $x_1 \approx - 4.780, \,\, x_2 = - 1.5$, and $x_3 = -0.471$, which confirms the prediction of the model. \\
Finally, for the equation
\b
\frac{8}{27} x^3 + \frac{4}{3} x^2 + 2 x + 1 = 0
\e
one has $c = 3b^2$ and, given that $a = (b/c)^3$, there must be a triple root given by $-3/c = -3/2$. This is so indeed.

\subsubsubsection{}
\vskip-0.73cm
\paragraph{\hskip1.5cm The Case of $\bm{a > 0, \,\, b > 0$}, and $\bm{c < 0}$ \\ \\}
$\hskip-.53cm$ This case is illustrated on Figure 4a. Consider first the sub-case of $c^2 > 4b$. The roots $\alpha_{1,2}$ in (\ref{alfa12}) have different signs ($\alpha_1 > 0$). There can be only one curve $\alpha_1 x^3 + 1$, with positive $\alpha_1$, tangent to the curve $-bx^2 - cx$. This occurs at point $\chi = (c/b) ( - 1 - \sqrt{1 - 3b/c^2}) < -c/b < 0$ which is the smaller root of equation (\ref{chi}). The other root, $(c/b)( - 1 + \sqrt{1 - 3b/c^2}) < -c/b > 0$, is associated with the irrelevant $\alpha_2 < 0$. \\
When $a$ is smaller than $\alpha_1$, the biggest root is positive and is between $\chi$ and $\sigma_1 = [-c/(2b)] (1 + \sqrt{1-4b/c^2})$ which is the bigger root of $ - b x^2 - c x - 1 = 0$. The middle root is also positive and is between $\sigma_2 = [-c/(2b)] (1 - \sqrt{1-4b/c^2})$ (which is the smaller root of $- b x^2 - c x - 1 = 0$) and $\chi$. The third root is negative and is smaller than $-\sqrt[3]{1/a}$. \\
If $a = \alpha_1$, there is a positive double root $\chi$ and a negative root smaller than $-\sqrt[3]{1/a}$. \\
Finally, if $a > \alpha_1$, then the cubic discriminant $\Delta_3$ is negative and the equation has a negative root smaller than $-\sqrt[3]{1/a}$ together with two complex roots. \\
When $c^2 = 4b$, the bigger root in (\ref{alfa12}) is $\alpha_1 = 0$. The other one is $\alpha_2 = c^3/54 < 0$. Thus the curve $a x^3 + 1$ with $a = \alpha_1 = 0$ is tangent to the curve $-bx^2 - cx$ at point $-c/(2b)$ where the maximum of $-bx^2 - cx$ occurs (the maximum in this case is $c^2/4b = 1$). As $a$ cannot be zero (one has to have an equation of degree 3), the cubic equation has a negative root smaller than $-\sqrt[3]{1/a}$ and two complex roots.\\
If $3b \le c^2 < 4b$, the roots $\alpha_{1,2}$ will have the same signs. The co-efficient in the term linear in $a$ is $2c(9b - 2c^2)$ and in the case of $3b \le c^2 < 4b$, given that $b$ is positive, its sign will depend on the sign of $c$ only. Thus, for negative $c$ the roots $\alpha_1$ and $\alpha_2$ will be negative and thus irrelevant.  No curve $a x^3 + 1$ with  $a > 0$ can intersect in the first quadrant the curve $-bx^2 - cx$ with $c^2 \le 4b$. The cubic discriminant $\Delta_3$ is non-negative. The roots of the cubic equation are as in the latter case: two complex and one negative and smaller than $-\sqrt[3]{1/a}$. \\
When $c^2 < 3b$, the discriminant $\Delta_3$ is negative and no curve $ax^3 + 1$ with any $a$ could be tangent to the curve $-bx^2 - cx$. The roots of the cubic equation are, again, two complex and one negative and smaller than $-\sqrt[3]{1/a}$. \\
This case will be illustrated with the following four examples. \\
Consider first the equation
\b
x^3 + x^2 - 5x +1 = 0.
\e
The roots (\ref{alfa12}) are $\alpha_1 = 205/27 + (44/27) \sqrt{22} \approx 15.236$ $\alpha_2 = 205/27 - (44/27) \sqrt{22} \approx -0.051$ (irrelevant). The relevant root $\chi = \chi_1$ is $5-\sqrt{22} \approx 0.310$. The roots of the equation $-bx^2 - cx - 1 = 0$ are $\sigma_1 = 5/2 +(1/2) \sqrt{21} \approx 4.791$ and $\sigma_2 = 5/2 - (1/2) \sqrt{21} \approx 0.209$. As the given $a$ is equal to 1 and smaller than $\alpha_1 \approx 15.236$, and as the given $b$ and $c$ are such that $c^2 > 4b$, then the roots of the equation are as follows: a positive root $x_1$ between $\chi \approx 0.310$ and $\sigma_1 \approx 4.791$, another positive root ($x_2$) between $\sigma_2 \approx 0.209$ and $\chi \approx 0.310$, and a negative root $x_3$ smaller than $-\sqrt[3]{1/a} = -1$. The roots of the equation are: $x_1 \approx 1.655$, $x_2 \approx 0.211$, and $x_3 \approx -2.866$ which agrees with the prediction. \\
Next, for the equation
\b
16x^3 + x^2 - 5x +1 = 0,
\e
the corresponding equation $\Delta_3 = 0$ has the same roots $\alpha_1 = 205/27 + (44/27) \sqrt{22} \approx 15.236$ and $\alpha_2 = 205/27 - (44/27) \sqrt{22} \approx -0.051$ (irrelevant). The relevant root $\chi = \chi_1$ is $5-\sqrt{22} \approx 0.310$ is the same. The roots of the equation $-bx^2 - cx - 1 = 0$ are also the same: $\sigma_1 = 5/2 +(1/2) \sqrt{21} \approx 4.791$ and $\sigma_2 = 5/2 - (1/2) \sqrt{21} \approx 0.209$. As the given $a$ is now equal to 16 and greater than $\alpha_1 \approx 15.236$, and as the given $b$ and $c$ are still such that $c^2 > 4b$, then the roots of the equation are as follows: a positive root $x_1$ between $\chi \approx 0.310$ and $\sigma_1 \approx 4.791$, another positive root ($x_2$) between $\sigma_2 \approx 0.209$ and $\chi \approx 0.310$, and a negative root $x_3$ smaller than $-\sqrt[3]{1/a} \approx -0.397$. The roots of the equation are: $x_{1,2} \approx 0.303 \pm 0.375\,i$ and $x_3 \approx 0.669$ which also agrees with the prediction.\\
As a further example, for the equation
\b
x^3 + 2x^2 - \frac{13}{5}x +1 = 0
\e
the corresponding equation $\Delta_3 = 0$ has roots $\alpha_1 = -1456/3375 + (38/3375) \sqrt{19} \approx -0.382$ and $\alpha_2 = - 1456/3375- (38/3375) \sqrt{19} \approx -0.480$. Both are irrelevant, since the only curves $ax^3 + 1$ that can be tangent to $-bx^2 - cx$ are those with $a < 0$ while the considered $a$ is positive. The equation should have two complex roots and a negative root smaller than $-\sqrt[3]{1/a} = -1$. This is the case indeed: $x_{1,2} \approx 0.492 \pm 0.305 \, i$ and $x_3 \approx - 2.984$. \\
Finally, the equation
\b
3x^3 + 2x^2 - x +1 = 0
\e
is with $c^2 < 3b$. The discriminant $\Delta_3$ is negative. Thus, there should be a root smaller than $-\sqrt[3]{1/a} \approx -0.693$ and two complex roots. This is the case indeed: $x_1 \approx -1.185$ and $x_{1,2} \approx 0.259 \pm 0.463 \, i$. \\

\subsubsubsection{}
\vskip-0.73cm
\paragraph{\hskip1.5cm The Case of $\bm{a > 0, \,\, b < 0$}, and $\bm{c >0}$ \\ \\}
$\hskip-.53cm$ Given that $b$ is negative, real roots (\ref{alfa12}) always exist and they always are with opposite signs. The relevant one is $\alpha_1 > 0$. The corresponding root of (\ref{chi}) is $\chi = \chi_1 =  (c/b) ( - 1 + \sqrt{1 - 3b/c^2}) > - c/b > 0$. The curve $ax^3 + 1$ with $a = 0$ intersects in the first quadrant the curve $-bx^2 - cx$ with $b < 0$ and $c >0$ at point $\sigma_1 = [-c/(2b)] (1 + \sqrt{1-4b/c^2})$ which is the bigger root of the equation $-bx^2 - cx - 1=0$. \\
Then the roots of the cubic equation are as follows. There is always one negative root between $-\sqrt[3]{1/a}$ and the origin. If $a > \alpha_1$, then the other two roots are complex. If $a = \alpha_1$, then in addition to the negative root between $-\sqrt[3]{1/a}$ and 0, there is a positive double root at $\chi$. If $a < \alpha_1$, the roots are: one negative root between $-\sqrt[3]{1/a}$ and 0; one positive root between $\sigma_1$ and $\chi$; and another positive root greater than $\chi$ --- see Figure 4b. \\
The equation
\b
x^3 - x^2 + x +1 = 0
\e
illustrates the case of $a > \alpha_1$: one has $a =1 > \alpha_1 = 5/27 \approx 0.185$. The method predicts one negative root between $-\sqrt[3]{1/a} = -1$ and 0 and two complex roots. The roots are $x_1 \approx -0.544$ and $x_{2,3} \approx 0.772 \pm 1.115\, i$. \\
The equation
\b
\frac{1}{9} x^3 - x^2 + x +1 = 0
\e
is in the category of $a = 1/9 < \alpha_1 = 5/27 \approx 0.185$. The locus of the corresponding double root is $\chi = \chi_1 = 3$. The bigger root of $-bx^2 - cx - 1 = 0$ is $\sigma_1 = 1/2 + (1/2) \sqrt{5} \approx 1.618$. The prediction of the method is for a negative root between $-\sqrt[3]{1/a} \approx -2.080$ and 0; a positive root between $\sigma_1 \approx 1.618$ and $\chi = 3$; and another positive root greater than $\chi = 3$. Indeed, the roots are: $x_1 \approx -0.607, \,\, x_2 \approx 1.932$, and $x_3 \approx 7.674$.

\begin{figure}[!ht]
\centering
\subfloat[\scriptsize ]
{\label{F4a}\includegraphics[height=6cm, width=0.48 \textwidth]{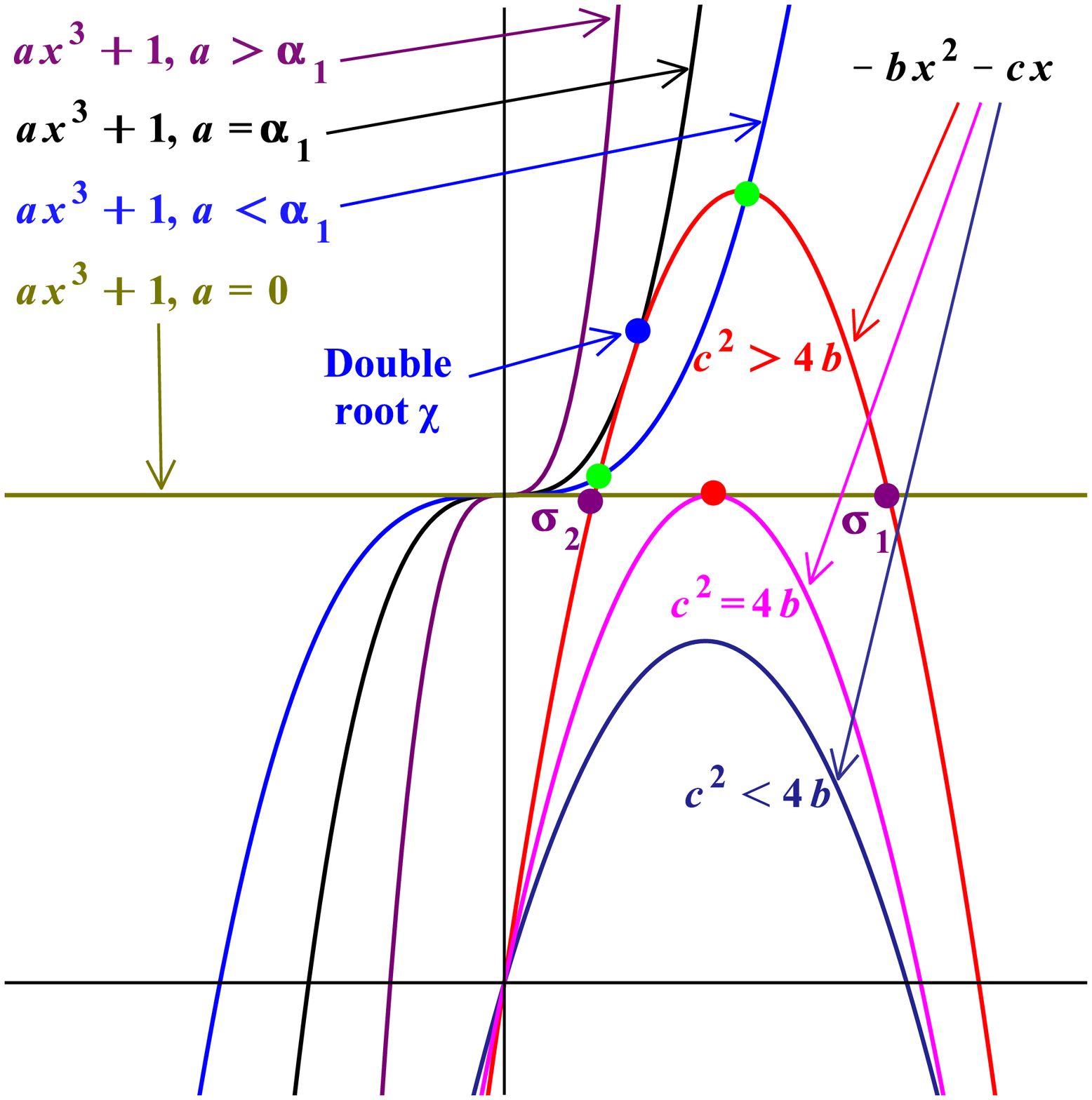}}
\quad
\subfloat[\scriptsize ]
{\label{F4b}\includegraphics[height=6cm,width=0.48 \textwidth]{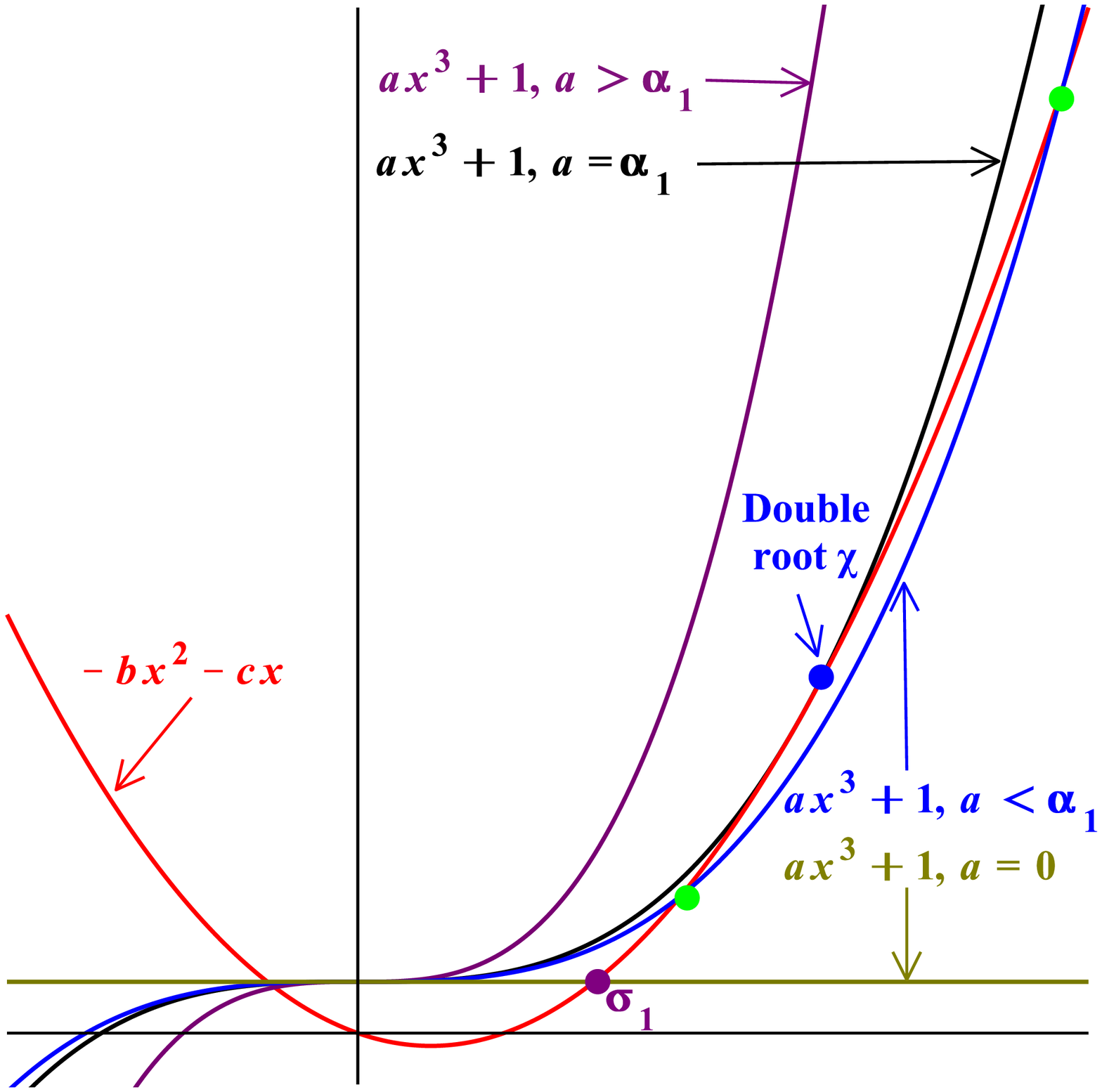}}
\end{figure}

\addtocounter{subfigure}{+2}

\begin{figure}[!ht]
\centering
\subfloat[\scriptsize ]
{\label{F4c}\includegraphics[height=6cm,width=0.54\textwidth]{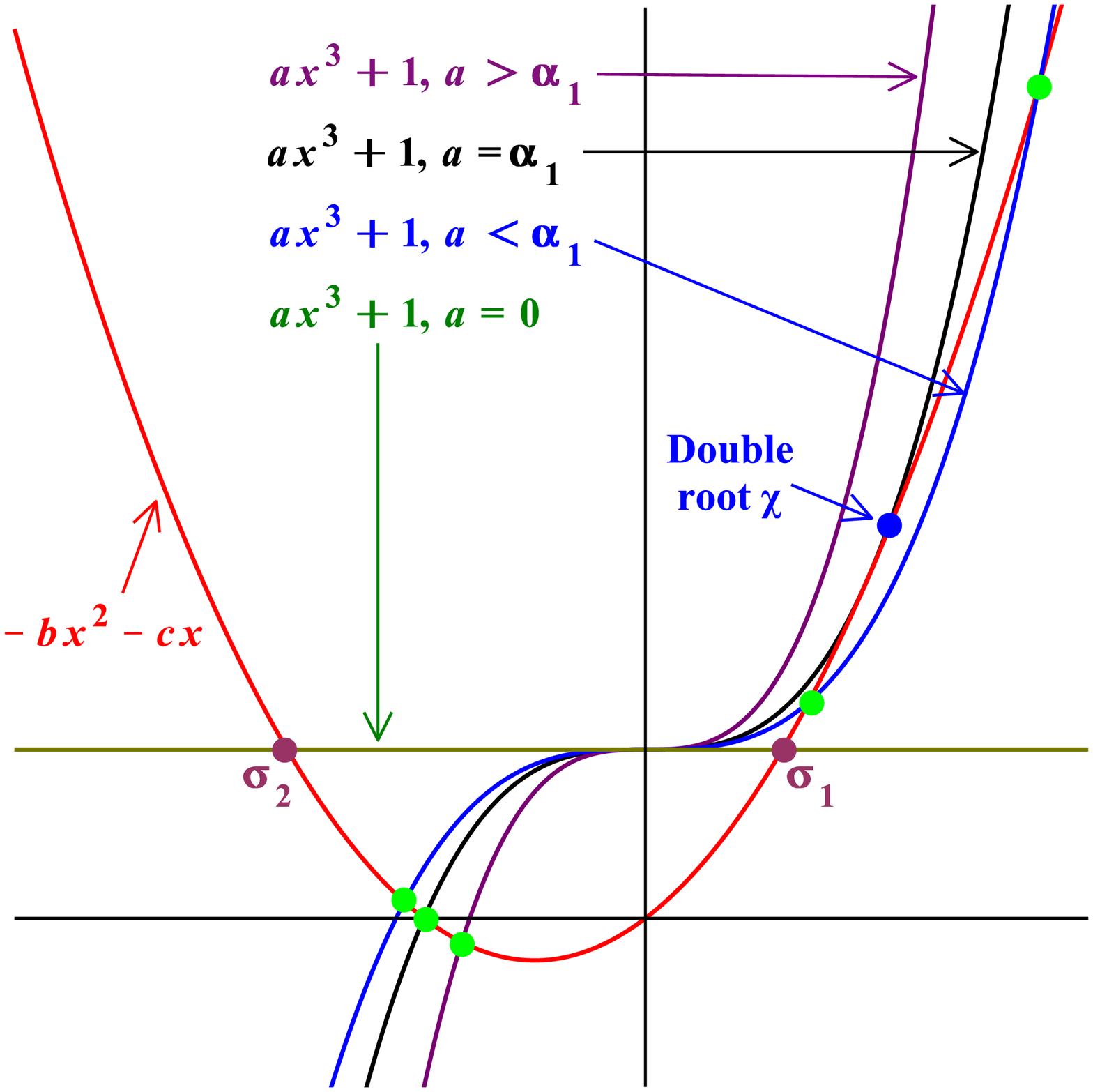}}
\caption{\footnotesize The cubic equation $ax^3 + 1 = -bx^2 - cx$ --- the cases of $a > 0, \,\, b > 0, \,\, c < 0$ (Figure 4a); $a > 0, \,\, b < 0, \,\, c > 0$ (Figure 4b); and $a > 0, \,\, b < 0, \,\, c < 0$ (Figure 4c). }
\label{Figure4}
\end{figure}

\subsubsubsection{}
\vskip-0.73cm
\paragraph{\hskip1.5cm The Case of $\bm{a > 0, \,\, b < 0$}, and $\bm{c < 0}$ \\ \\}
$\hskip-.53cm$ Given that $b$ is again negative, the roots (\ref{alfa12}) are again always real and always with opposite signs. The relevant one is $\alpha_1 > 0$. The corresponding root of (\ref{chi}) is $\chi = \chi_1 =  (c/b) ( - 1 + \sqrt{1 - 3b/c^2}) > - c/b > 0$. The curve $ax^3 + 1$ with $a = 0$ intersects in the first quadrant the curve $-bx^2 - cx$ with $b < 0$ and $c >0$ at point $\sigma_1 = [-c/(2b)] (1 + \sqrt{1-4b/c^2})$ which is the bigger root of the equation $-bx^2 - cx - 1=0$. \\
Then the roots of the cubic equation are as follows. There is always one negative root between $\mbox{min} \left( -\sqrt[3]{1/a}, -c/b \right)$ and $\mbox{max} \left( -\sqrt[3]{1/a}, -c/b \right)$. If $a > \alpha_1$, then the other two roots are complex. If $a = \alpha_1$, then, in addition to the negative root, there is a positive double root at $\chi$. If $a < \alpha_1$, the roots are: one negative root between $\mbox{min} ( -\sqrt[3]{1/a}, -c/b )$ and $\mbox{max} ( -\sqrt[3]{1/a}, -c/b )$; one positive root between $\sigma_1$ and $\chi$; and another positive root greater than $\chi$ --- See Figure 4c. \\
As an example, consider the equation
\b
2x^3 - x^2 - x + 1 = 0.
\e
The roots (\ref{alfa12}) are $\alpha_2 = -5/27$ (irrelevant) and $\alpha_1 = 1$. The corresponding loci of the double roots are $\chi_2 = -3$ (irrelevant) and $\chi = \chi_1 = 1$. The bigger root of $-bx^2 - cx - 1 = 0$ is $\sigma_1 = (1/2)(-1+\sqrt{5}) \approx 0.618$. Also $-\sqrt[3]{1/a} \approx -0.794 > -c/b = -1$. The given $a = 2$ is greater than $\alpha_1 = 1$, thus the equation must have one negative root between $c/b = -1$ and $-\sqrt[3]{1/a} \approx -0.794$ and two complex roots. This is so indeed: the roots are $x_1 \approx -0.829$ and $x_{2,3} \approx 0.665 \pm 0.401 \, i$. \\
The equation
\b
\frac{1}{2} x^3 - x^2 - x + 1 = 0
\e
is another example chosen so that one again has: $\alpha_2 = -5/27$ (irrelevant), $\alpha_1 = 1, \,\,  \chi_2 = -3$ (irrelevant), $\chi = \chi_1 = 1, \,\, \sigma_1 = (1/2)(-1+\sqrt{5}) \approx 0.618$. This time $-\sqrt[3]{1/a} \approx -1.260 < -c/b = -1$.
The given $a = 1/2$ is now smaller than $\alpha_1 = 1$, thus the equation must have one negative root between  $-\sqrt[3]{1/a} \approx -1.260$ and $c/b = -1$ and two positive roots --- one between $\sigma_1 \approx 0.618$ and $\chi = 1$ and another one greater than $\chi = 1$. The roots are: $x_1 \approx -1.170, \,\, x_2 \approx 0.689$, and $x_3 \approx 2.481$ --- in their predicted bounds.

\subsubsubsection{}
\vskip-0.73cm
\paragraph{\hskip1.5cm The Four Cases with $\bm{a < 0}$ \\ \\}
$\hskip-.53cm$ The analysis of these four cases is completely analogous as there is symmetry (reflection with respect to the ordinate) between them and the four cases already studied (one only needs to replace $c$ by $-c$ when $a$ is replaced by $-a$). \\
That is, the case of $a < 0, \,\, b > 0$, and $c > 0$ is analogous to the case of $a > 0, \,\, b > 0$, and $c < 0$ (Figure 4a); the case of $a < 0, \,\, b > 0$, and $c < 0$ is a complicated case analogous to the case of $a > 0, \,\, b > 0$, and $c > 0$ (Figure 3); the case of $a < 0, \,\, b < 0$, and $c > 0$ is analogous to the case of $a > 0, \,\, b < 0$, and $c < 0$ (Figure 4c); and, finally, the case of $a < 0, \,\, b < 0$, and $c < 0$ is analogous to the case of $a > 0, \,\, b < 0$, and $c > 0$ (Figure 4b).

\subsection{An Example with Equation of Degree 5}

\begin{figure}[!ht]
\centering
{\label{F5}\includegraphics[height=8cm, width=0.58\textwidth]{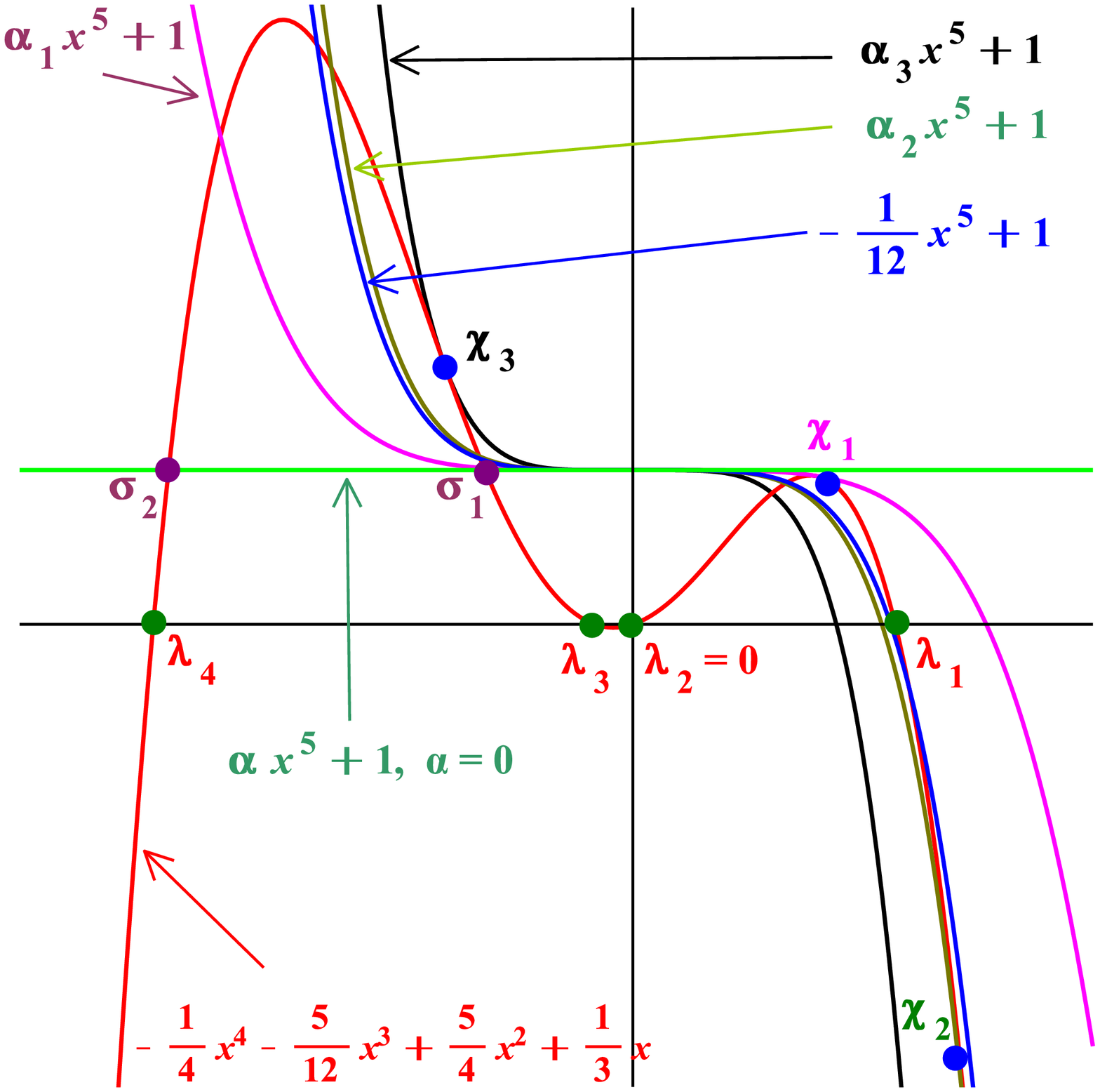}}
\caption{\footnotesize{The quintic equation $ - \frac{1}{12} x^5 + \frac{1}{4} x^4 + \frac{5}{12} x^3 - \frac{5}{4} x^2 - \frac{1}{3}x + 1 = 0$. }}
\label{Figure5}
\end{figure}

Consider the quintic equation
\b
- \frac{1}{12} \, x^5 + \frac{1}{4} \, x^4 + \frac{5}{12} \, x^3 - \frac{5}{4} \, x^2 - \frac{1}{3} \, x + 1 = 0
\e
and split as follows
\b
- \frac{1}{12} \, x^5  + 1 \,\, = \,\,  - \frac{1}{4} \, x^4 - \frac{5}{12} \, x^3 + \frac{5}{4} \, x^2 + \frac{1}{3} \, x.
\e
Setting the discriminant of the quintic
\b
\alpha \, x^5 + \frac{1}{4} \, x^4 + \frac{5}{12} \, x^3 - \frac{5}{4} \, x^2 - \frac{1}{3} \, x + 1
\e
equal to zero, results in an equation of degree four,
\b
3125 \, \alpha^4 + \frac{978851}{972} \, \alpha^3 + \frac{10519165}{186624} \, \alpha^2 - \frac{1149707}{1119744} \, \alpha - \frac{76507}{2985984} = 0,
\e
the roots of which are:  $\alpha_1 \approx -0.015, \,\, \alpha_2 \approx -0.089, \,\, \alpha_3 \approx -0.24$, and $\alpha_4 \approx  0.024$. As the given $a$ is $-1/12 \approx -0.083 < 0$, the root $\alpha_4 > 0$ is irrelevant. One also has $\alpha_3 < \alpha_2 < a < \alpha_1$.
The double roots $\chi_j$, associated with $\alpha_j$, are the roots of the equation
\b
\frac{1}{4} \, x^4 + 2 \times \frac{5}{12} \, x^3 + 3 \times \left(- \frac{5}{4} \right) \, x^2 + 4 \times \left(- \frac{1}{3} \right) \, x + 5 \times 1 = 0,
\e
namely $\chi_1 \approx 1.192, \,\, \chi_2 \approx 2.398, \,\, \chi_3 \approx - 1.229,$ and the irrelevant $\chi_4$ is $-5.695$ approximately. \\
Next, solve the equations
\b
\frac{1}{4} \, x^4 + \frac{5}{12} \, x^3 - \frac{5}{4} \, x^2 - \frac{1}{3} \, x + 1 = 0
\e
to determine the points at which the curve $ax^5 + 1$ with $a = 0$ intersects the curve $- (1/4) x^4 - (5/12) x^3 + (5/4) x^2 + (1/3) x$. The solutions are: $\sigma_1 \approx -0.965, \,\, \sigma_2 \approx -3.028$, and $\sigma_{3,4} \approx 1.163 \pm 0.124 \, i$. \\
The prediction of the model is for one negative root between $\sigma_2 \approx -3.028$ and  $\chi_3 \approx -1.229$; one negative root between $\chi_3 \approx -1.229$  and   $\sigma_1 \approx -0.965$; one positive root between  $0$   and   $\chi_1 \approx 1.192$; one positive root between  $\chi_1 \approx 1.192$  and   $\chi_2 \approx 2.398$; one positive root bigger than $\chi_2 \approx 2.398$ --- see Figure (\ref{Figure5}). \\
The roots are: $x_1 = -2, \,\, x_2 = -1, \,\, x_3 = 1, \,\, x_4 = 2, \,\, x_5 = 3$ --- within the predicted bounds.

\subsection{Recursive Application of the Method. An Example with Equation of Degree 7}

Consider the equation of seventh degree
\b
\frac{16}{3} \, x^7  - \frac{52}{3} \, x^6 + \frac{14}{3} \, x^5 + \frac{77}{3} \, x^4 - \frac{77}{6} \, x^3 - \frac{28}{3} \, x^2 + \frac{17}{6} \, x + 1 = 0.
\e
and see Figure 6. \\
In order to find the number of positive and negative roots of this equation, together with their bounds, one cannot consider setting its discriminant equal to zero as the resulting equation for the unknown $\alpha$ (which replaces the co-efficient $16/3$ of $x^7$) will be of degree 6 and not solvable analytically. One has to  proceed by applying the method twice. Firstly, re-write the given equation $p_7(x) = 0$ as $\lambda(x) - \mu(x) = 0$ with
\b
\lambda(x) & \!\!\! = \!\!\! & \frac{16}{3} \, x^7 + 1, \\
\mu(x) & \!\!\! = \!\!\! & - \frac{17}{6} \, x \left( -\frac{104}{17} \, x^5 + \frac{28}{17} \, x^4 + \frac{154}{17} \, x^3 - \frac{77}{17} \, x^2
- \frac{56}{17} \, x + 1 \right) \equiv - \frac{17}{6} \, x \, \hat{\mu}(x).
\e
The only real root of $\lambda(x) = 0$ is $\lambda_1 = - \sqrt[7]{3/16} \approx -0.787$. The roots of $\mu(x) = 0$ are those of $\hat{\mu}(x) = 0$ and zero. If the number of positive and negative ones among those can be determined, together with their bounds, then they can be used to attempt to determine the number of positive and negative roots (and their respective bounds) of the original equation. \\
Consider next the resulting equation of fifth degree $\hat{\mu} (x) = 0$, that is
\b
- \frac{104}{17} \, x^5 + \frac{28}{17} \, x^4 + \frac{154}{17} \, x^3 - \frac{77}{17} \, x^2 - \frac{56}{17} \, x + 1 = 0,
\e
and split further: $\hat{\mu}(x) = \rho(x) - \nu(x) = 0$, where:
\b
\rho(x) & \!\!\! = \!\!\! &  -\frac{104}{17} \, x^5 + 1, \\
\nu(x) & \!\!\! = \!\!\! & - \frac{28}{17} \, x^4 - \frac{154}{17} \, x^3 + \frac{77}{17} \, x^2 + \frac{56}{17} \, x.
\e
The only real root of $\rho(x) = 0$ is $\rho_1 = - \sqrt[5]{17/104} \approx 0.696$. The equation $\nu(x) = 0$ is of degree 4 and its roots are:
$\nu_1 \approx -5.908$, which is not shown on Figure 6b to avoid scaling of the graph, $\nu_2 \approx -0.412, \,\, \nu_3 = 0,$ and $\nu_4 \approx 0.821$. \\
Considering the quintic discriminant of $ \alpha \, x^5  + (28/17) \, x^4 + (154/17) \, x^3 - (77/17) \, x^2 - (56/17) \, x + 1$ and setting it to zero, results in an equation of degree four for the unknown $\alpha$ the roots of which are $\alpha_1 \approx -3.955, \,\, \alpha_2 \approx -10.021$ and the irrelevant ones $\alpha_3 \approx 0.072$ and $\alpha_4 \approx 60.743$ (given that the co-efficient of $x^5$ in $\rho(x)$ is $-104/17 \approx -6.118 < 0$ and $\alpha_{3,4}$ are positive). \\
The corresponding loci of the double roots of the equations $\alpha_j \, x^5  + (28/17) \, x^4$ \linebreak $ + \, (154/17) \, x^3 - (77/17) \, x^2 - (56/17) \, x + 1 = 0, \,\,\, j = 1,2,3,4, \,\,$ coincide with the roots of the equation
\b
\frac{2}{17} \, x^4 + 2 \times \, \frac{154}{17} \, x^3 + 3 \times \left( - \frac{77}{17} \right) \, x^2 + 4 \times \left( - \frac{56}{17} \right) \, x + 5 \times 1 = 0
\e
and are: $\chi_1 \approx 1.081, \,\, \chi_2 \approx - 0.753$, and the irrelevant $\chi_3 \approx -11.648$ and $\chi_4 \approx 0.320$. \\
The intersection points of the curve $\alpha x^5 + 1$ for which $\alpha = 0$ with the curve $\nu(x)$ are the roots of the equation $(28/17) \, x^4 + (154/17) \, x^3 - (77/17) \, x^2 - (56/17) \, x + 1 = 0$, namely: $\sigma_1 \approx -5.905, \,\, \sigma_2 \approx -0.559, \,\, \sigma_3 \approx 0.261,$ and $\sigma_4 \approx 0.704$.

\begin{figure}[!ht]
\centering
\subfloat[\scriptsize The original equation of seventh degree is split as $\lambda(x) - \mu(x) = 0$ with $\lambda(x) = \frac{16}{3} x^7 + 1$ and $\mu(x) = - \frac{17}{6} x \left( -\frac{104}{17} x^5 + \frac{28}{17} x^4 + \frac{154}{17} x^3 - \frac{77}{17} x^2 - \frac{56}{17} x + 1 \right)$. An equation of fifth degree, $\hat{\mu}(x) =  -\frac{104}{17} x^5 + \frac{28}{17} x^4$ $+ \frac{154}{17} x^3 - \frac{77}{17} x^2 - \frac{56}{17} x + 1 =0$, results. ]
{\label{F6a}\includegraphics[height=6cm, width=0.48\textwidth]{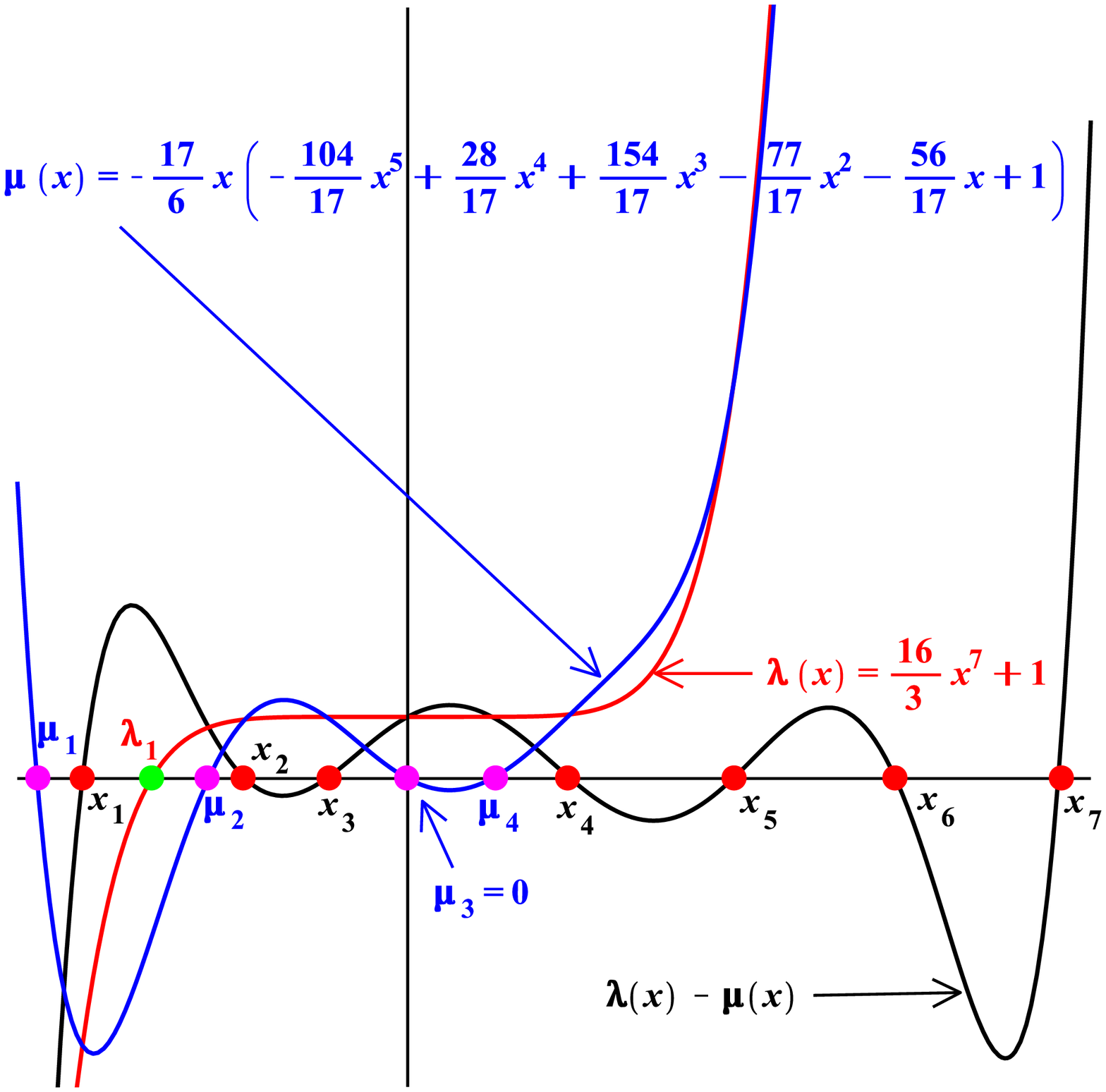}}
\quad
\subfloat[\scriptsize To find the roots of the resulting equation of fifth degree, $\hat{\mu}(x) =  -\frac{104}{17} x^5 + \frac{28}{17} x^4 + \frac{154}{17} x^3 - \frac{77}{17} x^2 - \frac{56}{17} x + 1 =0$, perform the further split $\hat{\mu}(x) = \rho(x) - \nu(x)$ with $\rho(x) = -\frac{104}{17} x^5 + 1 $ and $\nu(x) = - \frac{28}{17} x^4 - \frac{154}{17} x^3 + \frac{77}{17} x^2 + \frac{56}{17} x$. ]
{\label{F6b}\includegraphics[height=6cm,width=0.48\textwidth]{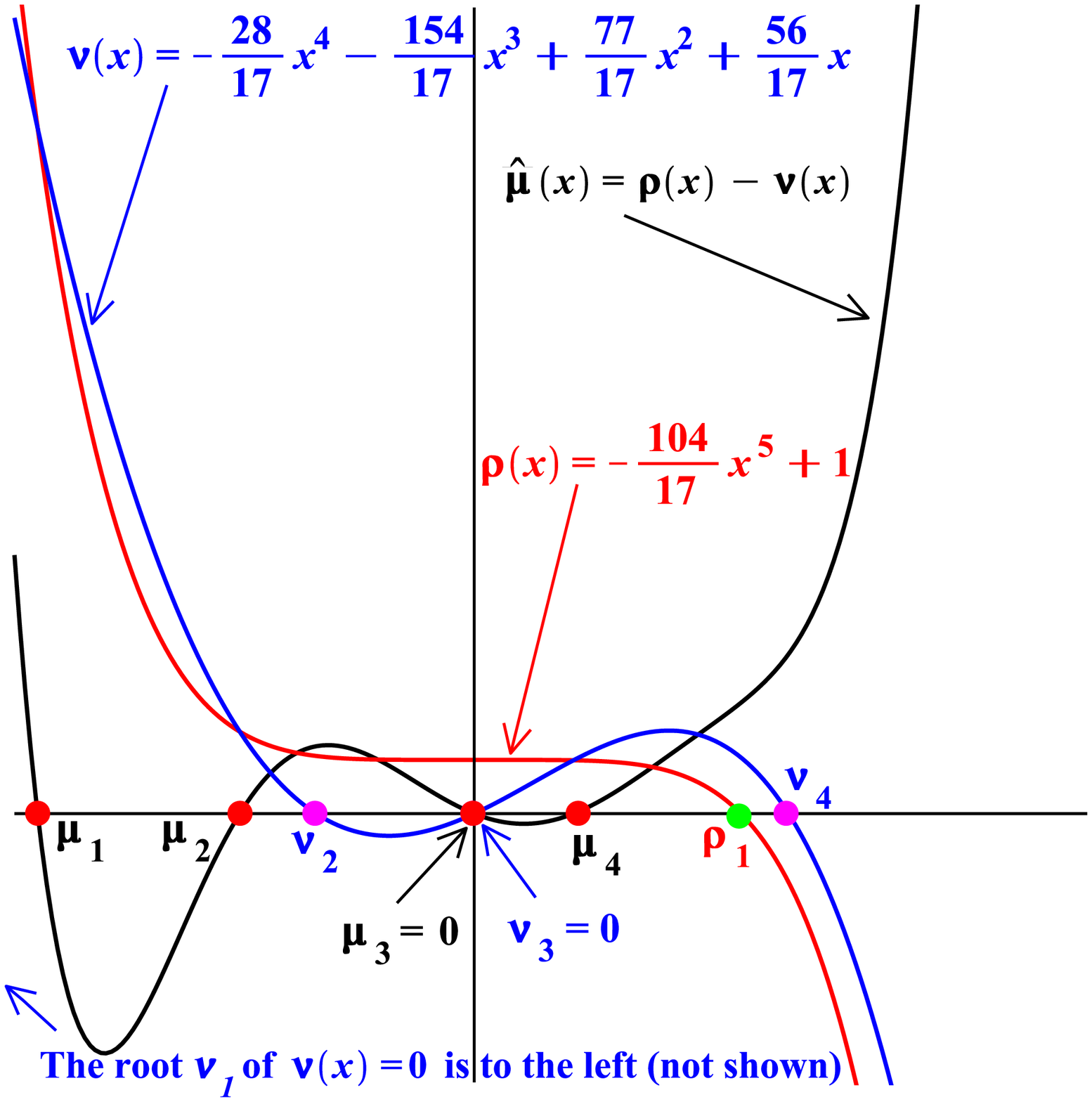}}
\label{Figure6}
\end{figure}

\addtocounter{subfigure}{+2}

\begin{figure}[!ht]
\centering
\subfloat[\scriptsize Since $\alpha_2 < - \frac{104}{17} < \alpha_1$, the roots of $\hat{\mu}(x) = \rho(x) - \nu(x) =  -\frac{104}{17} x^5 + \frac{28}{17} x^4 + \frac{154}{17} x^3 - \frac{77}{17} x^2 - \frac{56}{17} x + 1 = 0$ are as follows: one negative root smaller than $\chi_2$; another negative root between $\chi_2$ and $\sigma_2$; a positive root smaller than $- (104/17)^{-1/7}$; and two complex roots.]
{\label{F6c}\includegraphics[height=8cm, width=0.58\textwidth]{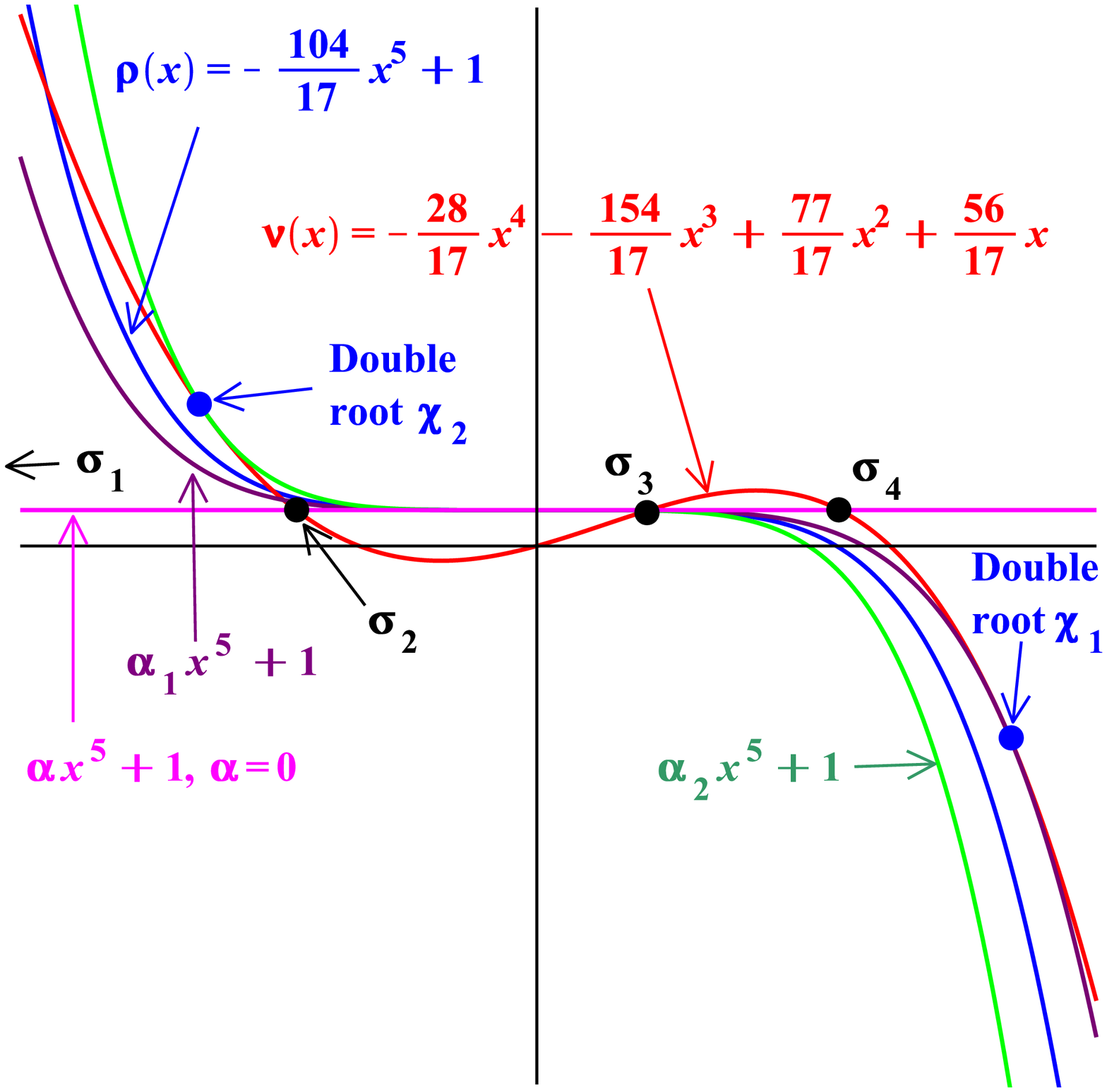}}
\caption{\footnotesize{Recursive application of the method for the equation $\frac{16}{3} x^7  - \frac{52}{3} x^6 + \frac{14}{3} x^5 + \frac{77}{3} x^4 - \frac{77}{6} x^3 - \frac{28}{3} x^2 + \frac{17}{6} x + 1 = 0$.}}
\label{Figure6_2}
\end{figure}

\n
The co-efficient of $x^5$ in $\rho(x)$ is $-104/17 \approx -6.118 < 0$ and this is greater than $\alpha_2 \approx -10.021$ and smaller than $\alpha_1 \approx -3.955$. Therefore, the roots of the equation $\mu(x) =  -[(17x)/6 ][-(104/17) x^5 + (28/17) x^4 + (154/17) x^3 - (77/17) x^2 - (56/17) x + 1] = 0$ are as follows: negative root $\mu_1$ smaller than $\chi_2 \approx - 0.753$; negative root $\mu_2$  between $\chi_2 \approx - 0.753$ and $\sigma_2 \approx -0.559$; $\,\, \mu_3 = 0$; positive root $\mu_4$ smaller than $\rho_1 = - (104/17)^{-1/7} \approx 0.547$; and two complex roots $\mu_{5,6}$. The actual roots are: $\mu_1 \approx -1.140, \,\, \mu_2 \approx -0.612, \,\, \mu_3 = 0, \,\, \mu_4 \approx 0.259$, and $\mu_{5,6} \approx 0.881 \pm 0.359\, i$  --- within their predicted bounds. \\
Returning with the roots of $\mu(x) = 0$ to the first split, $p_7(x) = \lambda(x) - \mu(x) = 0$, one can determine the following (see Figure 6a). The given equation $p_7(x) = 0$ has one negative root $x_1$ between $\mu_1$ and $\lambda_1 = -\sqrt[7]{3/16} \approx -0.787$. But given that $\mu_1  < \chi_2 \approx -0.753$, then all that can be said about the root $x_1$ is that it must be smaller than $\lambda_1 \approx -0.787$. The actual root of the equation is $x_1 = -1$. There can be no roots between $\lambda_1 \approx -0.787$ and $\mu_2$ as the function $\mu(x)$ is negative, while the function $\lambda(x)$ is positive there. Given that $\mu_2$ is between $\chi_2 \approx - 0.753$ and $\sigma_2 \approx -0.559$, then there could be no roots between $\lambda_1 \approx -0.787$ and  $\chi_2 \approx - 0.753$.  Next, consider the sub-interval between $\mu_2$ and 0. As $\mu_2$ is between $\chi_2 \approx - 0.753$ and $\sigma_2 \approx -0.559$ and given that $\mu(x)$ is of degree 5 and thus it can have up to four extremal points, then there could be either zero or two negative roots of the original equation between and $\sigma_2 \approx -0.559$ and the origin. The actual roots of the equation there are two: $x_2 = -1/2$ and $x_3 = -1/4$. Given that between $\mu_3 = 0$ and $\mu_4$, the function $\mu(x)$ is negative, while the function $\lambda(x)$ is positive, then there could be no intersection of these two curves in this interval and the equation $p_7(x) = 0$ cannot have roots in it. However,  given that $\mu_4$ is between 0 and $\rho_1 = - (104/17)^{-1/7} \approx 0.547$, what can be said about the positive roots of this equation is that there is either two, or four of them. The actual roots are $x_4 = 1/2, \,\, x_5 = 1, \,\, x_6 = 3/2,$ and $x_7 = 2$.

\section{Using the Split (\ref{p1}). An Example with Equation of Degree 9}

\begin{figure}[!ht]
\centering
{\label{F7}\includegraphics[height=8cm, width=0.58\textwidth]{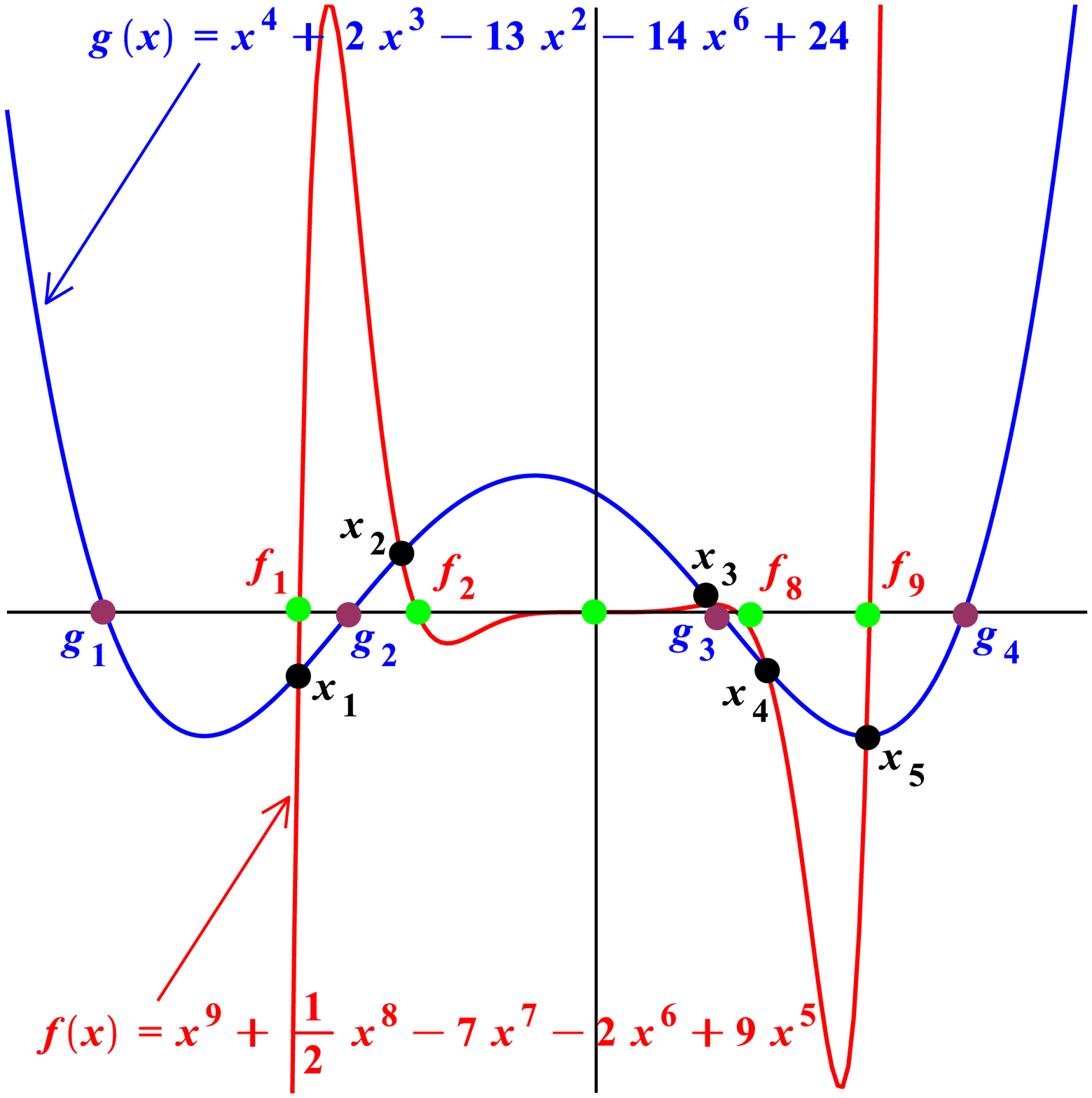}}
\caption{\footnotesize{The equation $x^9 + \frac{1}{2} x^8 - 7 x^7 - 2 x^6 + 9 x^5 - x^4 - 2 x^3 + 13 x^2 + 14 x - 24 = 0$. }}
\label{Figure7}
\end{figure}

Consider the following example as an equation of degree above 5 and up to and including 9:
\b
x^9 + \frac{1}{2} x^8 - 7 x^7 - 2 x^6 + 9 x^5 - x^4 - 2 x^3 + 13 x^2 + 14 x - 24 = 0.
\e
Using the split (\ref{p1}), this equation can be written as $f(x) - g(x) = 0$ with
\b
f(x) & = & \left( x^4 + \frac{1}{2} x^3 - 7 x^2 - 2 x + 9 \right) x^5, \\
g(x) & = & x^4 + 2 x^3 - 13 x^2 - 14 x + 24.
\e
The roots of the two equations $f(x) = 0$ and $g(x) = 0$ can be determined analytically. \\
For $f(x) = 0$, these are: $f_1 \approx -2.416, \,\, f_2 \approx -1.458, \,\, f_{3,4,5,6,7} = 0$ (zero is a quintuple root), $f_8 \approx 1.145$, and $f_9 \approx 2.230$. The origin is a saddle for $f(x)$. The first four derivatives of $f(x)$ at 0 are zero, while the fifth one is positive. Thus, $f(x)$ enters through the origin into the first quadrant from the third. The function $f(x)$ has two negative roots and two positive ones. Further, when $x \to -\infty, \,\, f(x) \to -\infty$, while when $x \to \infty, \,\, f(x) \to \infty$. It is essential to note that $f(x)$ can have up to 4 non-zero extremal points. Namely, $f'(x) = 9 x^8 + 4 x^7 - 49 x^6 - 12 x^5 + 45 x^4$ and setting this to zero gives the following extremal points:  a quadruple 0, together with the points: $-2.179, \,\, -1.210, \,\, 0.952$, and $1.992$. \\
The roots of $g(x) = 0$ are: $g_1 = -4, \,\, g_2 = -2, \,\, g_3 = 1,$ and $g_4 = 3$. At zero, one has $g(0) = 24 > 0$. The function $g(x)$ tends to $+\infty$ when $x \to \pm \infty$. It also has two positive roots and two negative ones. The function $g(x)$ can have up to three extremal points.  These are the roots of the equation $g'(x) = 4 x^3 + 6 x^2 - 26 x - 14 = 0$, namely the points $-3.193, \,\, -0.5,$ and $2.193$.\\
This very simple analysis allows the graphs of $f(x)$ and $g(x)$ to be easily sketched --- see Figure (\ref{Figure7}) --- and from the graph one can infer the following for the roots $x_i$ of the given equation. There can be no roots smaller than $g_1 = -4$ as in that region $f(x)$ and $g(x)$ have different signs and thus cannot intersect. There is one negative root between $g_1 = - 4$ and $f_1 \approx -2.416$. Given the number of extremal points of $f(x)$, it is not possible to have more negative roots in this sub-interval. There can be no roots between $f_1 \approx -2.416$ and $g_2 = 2$. There is one negative root between $g_2 = - 2$ and $f_2 \approx -1.458$. Again, due to the number of extremal points of $f(x)$, there can be no other negative roots in this sub-interval. There are no roots between $f_2 \approx -1.458$ and zero. There is one and no more positive roots between 0 and $g_3 = 1$. There can be no roots between $g_3 = 1$ and $f_8 \approx 1.145$. Between $f_8 \approx 1.145$ and $f_9 \approx 2.230$, there can be either zero or 2 positive roots. There can be no roots between $f_9 \approx 2.230$ and $g_4 = 3$. As $f(x)$ grows faster than $g(x)$, there can be no roots greater than $g_4 = 3$ either. The biggest positive root is therefore smaller than $f_9 \approx 2.230$. All roots are locked between $g_1 = -4$ and $f_9 \approx 2.230$.

\n
The actual roots of the equation are: $x_1 \approx - 2.426, \,\, x_2 \approx -1.591, \,\, x_3 \approx 0.948, \,\, x_4 \approx 1.388, \,\,  x_5 \approx 2.203, \,\, x_{6,7} \approx -0.917 \pm 0.837\, i$ and $x_{8,9} \approx 0.406 \pm 1.107\, i$. \\
For comparison, the Descartes rule of signs provides for either 1, or 3, or 5 positive roots and for either 0, or 2, or 4 negative roots. \\
The Lagrange bound provides that all roots are between $\pm $ max $ \left(1, \sum_{i=0}^{n-1} \vert a_i / a_n \vert \right) = \pm 72.5$.
Cauchy's theorem provides a stricter bound on all roots: they are locked between $\pm \Bigl( 1 + \Bigr.$ \raisebox{-1.2ex}{$\stackrel{\mbox{max}}{\scriptscriptstyle 0 \le k \le n-1}$} $\Bigl. \vert a_k \vert \Bigr) = \pm 25.$ \\

\section{A Different Perspective on the Split (\ref{p2})}

Given the equation $a_n \, x^n + a_{n-1} \, x^{n-1} + a_{n-2} \, x^{n-2} + \ldots + a_1 \, x + a_0 = 0$ which is of degree $n$ and thus the coefficient $a_n$ cannot be zero, one can instead consider an equation with $a_n$ set equal to 1:
\b
x^n + a_{n-1} x^{n-1} + a_{n-2} x^{n-2} + \ldots + a_1 x + a_0 = 0,
\e
which has the same roots.  \\
Then the split
\b
x^n + a_0 = - a_{n-1} x^{n-1} - a_{n-2} x^{n-2} - \ldots - a_1 x
\e
would allow the ``propagation" of the curve of fixed shape $x^n$ vertically until one finds the tangent points of $x^n + a_0$ with the curve $- a_{n-1} x^{n-1} - a_{n-2} x^{n-2} - \ldots - a_1 x$. In this case, one will again have to solve equations of degree one smaller than that of the given equation as $n-1$ is the highest power of $a_0$ in the discriminant $\Delta_n$ of $x^n + a_{n-1} x^{n-1} + a_{n-2} x^{n-2} + \cdots + a_1 x + a_0$. Then, in order to determine the number of positive and negative roots of the given equation and to also find their bounds, the given $y$-intercept $a_0$ will have to be compared to the real roots $\gamma_i$ of the equation $\Delta_n(\gamma) = 0$ in which $a_0$ has been replaced by $\gamma$ and treated as an unknown, while all other $a_j \,\, (j = 1,2, \ldots, n-1)$ are as given. This is an equation of degree $n-1$ and it is solvable analytically for $n \le 5$. For $n \ge 6$, this split should be applied recursively. The idea is very similar to the one studied in detail in this paper.

\end{document}